\documentclass[pdflatex,sn-mathphys-num]{sn-jnl}

\usepackage{graphicx}
\usepackage{multirow}%
\usepackage{amsmath,amssymb,amsfonts,soul}
\usepackage{amsthm}%
\usepackage{enumitem}
\usepackage{mathtools}
\usepackage{mathrsfs}%
\usepackage[title]{appendix}%
\usepackage{xcolor}%
\usepackage{textcomp}%
\usepackage{manyfoot}%
\usepackage{booktabs}%
\usepackage{algorithm}%
\usepackage{algorithmicx}%
\usepackage{algpseudocode}
\usepackage{listings} 
\usepackage{subfig}
\usepackage{fixltx2e}
\usepackage{bm}
\usepackage{accents}
\usepackage{bbm}
\usepackage{stmaryrd}
\usepackage{hyperref}
\usepackage{float}
\usepackage{titlesec}

\titleformat{\paragraph}
{\normalfont\normalsize\bfseries}{\theparagraph}{1em}{}
\titlespacing*{\paragraph}
{0pt}{3.25ex plus 1ex minus .2ex}{1.5ex plus .2ex}

\newcommand*{\dt}[1]{
\accentset{\mbox{\large\bfseries .}}{#1}}
\theoremstyle{thmstyleone}%
\theoremstyle{thmstyletwo}%

\theoremstyle{thmstylethree}%

\begin{document}

\title{Wavelet $s$-Wasserstein distances for $0<s\leqslant1$}


\author{Katy Craig}\email{kcraig@ucsb.edu}

\author{Haoqing Yu}\email{haoqingyu@ucsb.edu}

\affil{\orgdiv{Department of Mathematics}, \orgname{University of California, Santa Barbara}}


\abstract{
Motivated by classical harmonic analysis results characterizing H\"older spaces in terms of the decay of their wavelet coefficients, we consider wavelet methods for computing $s$-Wasserstein type distances. Previous work by Sheory (n\'e Shirdhonkar) and Jacobs showed that, for $0<s\leqslant1$, the $s$-Wasserstein distance $W_s$ between certain probability measures on Euclidean space is equivalent to a weighted $\ell^1$ difference of their wavelet coefficients. We demonstrate that the original statement of this equivalence is incorrect in a few aspects and, furthermore, fails to capture key properties of the $W_s$ distance, such as its behavior under translations of probability measures. Inspired by this, we consider a variant of the previous wavelet distance formula for which equivalence (up to an arbitrarily small error) does hold for $0<s<1$. We analyze the properties of this distance, one of which is that it provides a natural embedding of the $s$-Wasserstein space into a linear space. We conclude with several numerical simulations. Even though our theoretical result merely ensures that the new wavelet $s$-Wasserstein distance is   equivalent to the classical $W_s$ distance (up to an error), our numerical simulations show that the new wavelet distance succeeds in capturing the behavior of the exact $W_s$ distance under translations and dilations of probability measures. 
}




\maketitle

\section{Introduction}\label{sec1}
\subsection{Motivation and main contributions}\label{subsec1}
\noindent Over the past twenty years, \emph{Wasserstein distances} on the space of probability measures $\mathcal{P}(X)$ have attracted significant attention, in both pure and applied mathematics, due to their unique ability to respect the geometry on the underlying space $X$. However, a key drawback is that computing a Wasserstein distance between two measures involves solving a computationally expensive optimization problem. This has sparked a significant amount of research focused on either developing efficient algorithms for Wasserstein distances or finding easily computable substitutes that resemble them. For an excellent survey, see Peyré \& Cuturi \cite{bib17}.

When the underlying space is Euclidean, $X = \mathbb{R}^n$, Sheorey (né Shirdhonkar) \& Jacobs \cite{bib1} brought in the tool of wavelet systems to obtain an alternative to $s$-Wasserstein distances $W_s$ with $0<s\leqslant1$; see Section \ref{sec2} for backgrounds on $W_s$ and wavelets. For a subfamily of probability measures with densities on $\mathbb{R}^n$, the authors defined the original wavelet $s$-Wasserstein distance with $C_0\geqslant0$ and $C_1>0$ as
$$(W_{\text{wav,ori}}^{(\phi,\psi,C_0,C_1)})_{s}(\mu,\nu) \vcentcolon  = \sum_{\mathbf{k} \in \mathbb{Z}^n} C_0|(p_{\mu}-p_{\nu})_{(0,\mathbf{k})}^{(a)}| +\sum_{\substack{\epsilon\in\alpha,\\ j\geqslant 0,\  \mathbf{k}\in\mathbb{Z}^n}}C_12^{-j(s+\frac{n}{2})}|(p_{\mu}-p_{\nu})_{(\epsilon, (j, \mathbf{k}))}^{(d)}|,$$
where $(p_{\mu}-p_{\nu})_{(0,\mathbf{k})}^{(a)}$ and $(p_{\mu}-p_{\nu})_{(\epsilon, (j, \mathbf{k}))}^{(d)}$ are the approximation and detail coefficients, under the wavelet system $(\phi,\psi)$, of the difference densities $p_{\mu}-p_{\nu}$ for two measures $\mu$ and $\nu$; in particular, $C_0 = 0$ and $C_1 = 1$ were advocated and set to be the default choices in their paper. 

Due to the fast algorithm available for computing the discrete wavelet transform, a finite approximation to this wavelet distance can be computed in linear time $O(n)$ with respect to the number of sample points $n$ of the measures, contrasted with the $O(n^3\text{log}(n))$ complexity of the exact network simplex solver for $W_s$ \cite{bib17}. Also, the expression for $(W_{\text{wav,ori}}^{(\phi,\psi,C_0,C_1)})_{s}$ is in a closed form that involves only wavelet coefficients. Consequently, there is no need to tune any parameters in the computation to determine appropriate regularization or stopping criteria, making it simple to implement compared with other approximate algorithms. The explicit formula also leads to simpler theoretical estimates.

The main theoretical result in \cite{bib1} was that, for any $C_0\geqslant0$ and $C_1>0$, $(W_{\text{wav,ori}}^{(\phi,\psi,C_0,C_1)})_{s}$ is equivalent to $W_s$ for $0<s<1$\footnote{The $s=1$ case is subtle, see Sections \ref{subsec3.1} and \ref{subsec3.2} for details.}. This was established on the basis that certain regular wavelet systems can characterize the $s$-Hölder constraint set in the Kantorovich-Rubinstein dual formulation of $W_s$ for $0<s<1$. The characterizations of $s$-Hölder spaces, and more general Besov spaces, by certain regular wavelet systems have been well-known to the harmonic analysis community since 90s, due to the work of Meyer \cite[Chapter 6.4 \& Chapter 6.10]{bib4}, which was related further back to the classical Littlewood-Paley like characterizations of Besov spaces, due to Peetre \cite{bib31}. For these reasons, the main contribution of Sheory \& Jacobs was to explicitly estimate the size of some constants in the existing characterizations (see section \ref{subsec3.1}),  from which they conclude the equivalence of $(W_{\text{wav,ori}}^{(\phi,\psi,C_0,C_1)})_{s}$ and $W_s$.

However, in the present work, we demonstrate via a counterexample that the equivalence claimed in the original \cite{bib1} is invalid and pinpoint the exact part in the proof, as they combined their estimations, that led to the error. We show that this mistake also incurs problems on the numerical side, and give examples and intuition based on the computation of $(W_{\text{wav,ori}}^{(\phi,\psi,C_0,C_1)})_{s}$ when $C_0=0$ and $C_1=1$; see Sections \ref{subsec3.2}.

Next, we propose two approaches to fix the problem. First, we consider a naive solution by requiring $C_0>0$ and the domain to be compact. While this does restore the theoretical equivalence of $(W_{\text{wav,ori}}^{(\phi,\psi,C_0,C_1) })_{s}$ with $W_s$ for $0<s<1$ on compact domains, it still exhibits unsatisfactory numerical performance in comparison with $W_s$. We believe this is due to the artificial choice of the combination of $C_0$ and $C_1$, as well as the incompatibility of this formulation of the wavelet distance with translations, which is not shared by the actual $W_s$; see Section \ref{subsec3.3}.

Inspired by this, we derive the following new reformulation of wavelet $s$-Wasserstein distance, $(W_{\text{wav}}^{(\psi,j_0)})_s$, that is again equivalent (up to an arbitrarily small error) to $W_s$ for $0<s<1$ on compact domains, but avoids the issues in the last paragraph and shows superior numerical performance. Our main theorem is as follows:\\

\noindent \textbf{Theorem.} Fix a compact set $K\subseteq \mathbb{R}^n$, $0<s<1$, and a compactly supported, continuously differentiable orthonormal wavelet system $(\phi,\psi)$. For any  probability measures $\mu$, $\nu$ on $K$ with square-integrable densities $p_{\mu}$, $p_{\nu}$, define the \textbf{new wavelet $\pmb{s}$-Wasserstein distance} $(W_{\text{wav}}^{(\psi,j_0)})_s(\mu,\nu)$ \textbf{with the lowest level} $\pmb{j_0\in\mathbb{Z}}$ by
\begin{equation}
(W_{\text{wav}}^{(\psi,j_0)})_s(\mu,\nu) \vcentcolon  = \sum_{\substack{\epsilon\in\alpha,\\ j\geqslant j_0,\  \mathbf{k}\in\mathbb{Z}^n}}2^{-j(s+\frac{n}{2})}|(p_{\mu}-p_{\nu})_{(\epsilon, (j, \mathbf{k}))}^{(d)}|. \notag
\end{equation}
Then, for any $\eta>0$, there exists sufficiently small $j_0\in\mathbb{Z}$, depending on $s$, the domain $K$, and the system $(\phi,\psi)$, such that, for all $\mu$ and $\nu$,
\begin{equation}
E_{\psi,s}W_s(\mu, \nu)-\eta\leqslant (W_{\text{wav}}^{(\psi,j_0)})_s(\mu, \nu) \leqslant \tilde{E}_{\psi,s,j_0}W_s(\mu, \nu)\notag ,
\end{equation}
where the constants $E_{\psi,s}$ and   $\tilde{E}_{\psi,s,j_0}$ depend only on the wavelet system $(\phi,\psi)$, $s$, and $j_0$.\\

\noindent See Theorem \hyperref[Theorem 3.13]{3.13} for a detailed statement of the above result, with precise formulas for relevant constants. We further prove the tightness of the lower bound, examine the possibilities of extensions to larger classes of probability measures, and discuss the dependence of $(W_{\text{wav}}^{(\psi,j_0)})_s$ on $j_0$; see Section \ref{subsec3.4}. 

Beyond the numerical benefits of $(W_{\text{wav}}^{(\psi,j_0)})_s$ when computing the distance between a fixed pair of probability measures, a significant motivation for the study of $(W_{\text{wav}}^{(\psi,j_0)})_s$ arises when one has a dataset of $N$ probability measures and seeks to apply machine learning techniques to cluster or classify the measures with the Wasserstein distance. Traditionally, $N(N-1)/2$ expensive computations are needed to obtain all pairwise optimal transport distances, which can be computationally prohibitive when $N$ is large --- for example, in particle physics applications of Cai et al. \cite{bib29}, $N = 10^6$. A powerful solution to this problem is to consider an embedding of the dataset into a linear space with an approximate (or equivalent) distance. In the linearized optimal transport framework developed by Wang et al. \cite{bib10} for the 2-Wasserstein distance, the number of expensive computations can be reduced to $O(N)$, the cost of embedding, after which the remaining clustering/classification can be done efficiently in the embedded space. 

In our setting, the map $$p_{\mu}(\mathbf{x})d\mathbf{x}\mapsto((p_{\mu})_{(\epsilon, (j,\mathbf{k}))}^{(d)})_{\substack{\epsilon\in\alpha,\\j\geqslant j_0,\mathbf{k}\in\mathbb{Z}^n}},$$ arising from our new  wavelet $s$-Wasserstein distance, provides a natural embedding of probability measures on a compact domain with square integrable densities into a linear space, with a distance equivalent to the traditional $W_s$, $0<s<1$, up to a small error. There have been previous works that interpret the original wavelet $s$-Wasserstein distance as an embedding, like Saabni \& Zwilling \cite{bib32}. We generalize this perspective to our new wavelet $s$-Wasserstein distance and explain some unique advantages compared to the original embedding shown in the linearized optimal transport framework \cite{bib10} and the popular $\dot{H}^{-1}$ embeddings in Peyre \cite{bib27} and Greengard et al. \cite{bib28}, which were based on Loeper \cite{bib33}. For example, our map preserves the sparsity of the measures in the data set; see Section \ref{subsec3.5}. 

Following our theoretical study of the wavelet $s$-Wasserstein distance and the associated linear embedding, we offer detailed procedures for numerically computing all three formulations of the wavelet distances and compare their numerical performance with the traditional $W_s$. We do not expect any wavelet distance to provide a method of exactly \emph{approximating} $W_s$ to a certain degree of accuracy, as the distances are merely equivalent. However, one of the key benefits of the linearized optimal transport framework, as pointed out in Moosmüller \& Cloninger \cite{bib12}, is that it precisely coincides with the $W_2$ distance on datasets generated by translations and dilations, which are fundamental types of measure transport. We use this property as the baseline to rate the behavior of three formulations of the wavelet distances and numerically show that only our new distance $(W_{\text{wav}}^{(\psi,j_0)})_s$ exhibits this type of similarity with $W_s$, after appropriate normalization; see Section \ref{subsec4.2}. We also conduct further numerical studies of how the performance of the new $(W_{\text{wav}}^{(\psi,j_0)})_s$ depends on the choice of wavelet system and the number of coefficients being summed; see Section \ref{subsec4.3}. Additionally, since $W_s\to W_1$ as $s\to1$, where $W_1$ is the well-known Earth Mover's Distance, one might still be interested in the performances of the wavelet distances when $s=1$, despite the fact this falls outside the scope of the theoretical results. We therefore include $s=1$ in all of our simulations in Section \ref{sec4}, and still observe good numerical results for the new $(W_{\text{wav}}^{(\psi,j_0)})_s$ when $s=1$.

\subsection{Related work}
On the theoretical level, in Leeb \& Coifman \cite{bib21}, Leeb \cite{bib24}, and Weed \& Berthet \cite{bib2}, the authors briefly mentioned the wavelet characterization of $s$-Hölder spaces and made connection with the Kantorovich Rubinstein dual formulation for $W_s$, $0<s<1$, where a subset of the $s$-Hölder space is the constraint set; this is, as explained earlier, the core idea in the original formulation of the wavelet $s$-Wasserstein distance in \cite{bib1}. As for final objectives, \cite{bib21} and \cite{bib24} are mainly concerned with a high-level analysis of certain Hölder–Lipschitz norms and their duals on spaces
with semigroups and with tree distances, respectively, while \cite{bib2} focused on minimax rates for nonparametric estimation in $W_s$ ($s\geqslant1$). Also somewhat related to the present line of research was the work by Weed \& Bach \cite{bib22}, in which a dyadic transport technique was used to bound $W_s$, with the goal of studying the convergence rate of estimators in $W_s$ ($s\geqslant1$).

Although there are some connections between the above-described papers and \cite{bib1}, it is clear that their primary objectives were different from ours: to the best of our knowledge, no other work has discovered the problem of the original wavelet $s$-Wasserstein distance and provided a resolution as we do.

Similarly, from a numerical view, there exists a large amount of research that used the formula  $(W_{\text{wav,ori}}^{(\phi,\psi,C_0,C_1)})_{s}$ with $C_0 = 0, C_1 = 1$ in computation as a surrogate for $W_s$; see Sandler \& Lindenbaum \cite{bib25}, Zelesko et al. \cite{bib23}, Lai et al. \cite{bib26}, etc. However, none of these works conducted a detailed study on the performance of $(W_{\text{wav,ori}}^{(\phi,\psi,C_0,C_1)})_{s}$ and compared with $W_s$. Even in the original work \cite{bib1}, the numerical simulations were done on a very small discrete domain that avoids the exposure of certain issues of $(W_{\text{wav,ori}}^{(\phi,\psi,C_0,C_1)})_{s}$. In this way, ours is the first work to conduct a thorough numerical analysis of qualitative properties of various wavelet $s$-Wassertsein distances and contrast the behavior with $W_s$.

From the perspective of linearization of Wasserstein distance, along with the embeddings previously mentioned, there are also a variety of embeddings for the 1-Wasserstein distance, due to early works by Indyk \& Thaper \cite{bib19} and Andoni et al. \cite{bib20}. Nevertheless, these embeddings only consider measures supported on a fixed, bounded grid; in view of the type of measures being embedded, their methods are more different from the ones we shall discuss and compare in this paper, where the embedding maps are well-defined at the continuum level and thus independent of discretizations.

\subsection{Future work}
There are several directions for future work. On the theoretical side, we leave  two open directions: the tightness of the upper bound in the equivalence of our major theorem (see Remark \hyperref[Remark 3.17]{3.17}) and the extent to which the equivalence can be extended to a larger family of probability measures (see Remark \hyperref[Remark 3.18]{3.18}).

On the numerical side, the computational performance of our new wavelet $s$-Wasserstein distance between measures in higher dimensions $n>1$ is left for future exploration. Moreover, it would also be interesting to apply the linear embedding introduced in Section \ref{subsec3.5} to clustering and classifying real datasets.

\section{Preliminaries}\label{sec2}
\subsection{Basic notation}\label{subsec2.1}
Bold letters denote points in $\mathbb{R}^n$ or $\mathbb{Z}^n$ and normal letters with subscripts indicate their components; i.e. $\mathbf{x} = (x_1,\cdots,x_n)$. $B_r$ means an open ball of radius $r>0$ in Euclidean spaces. The notation $\llbracket n \rrbracket$ represents the set of integers from 1 to $n$.

Given $s>0$ and $E\subseteq \mathbb{R}^n$ closed, let $\mathcal{P}(E)$ be the set of Borel probability measures on $E$, $\mathcal{P}_s(E)$ be the set of elements in $\mathcal{P}(E)$ with finite $s$-moment, $M_s(\mu) =\int_{E}|\mathbf{x}|^sd\mu(\mathbf{x}) 
< +\infty$, and $\mathcal{P}_{c,\text{ab2}}(E)$ be the set of elements in $\mathcal{P}(E)$ that are compactly supported and absolutely continuous to the $n$-dimensional Lebesgue measure with square-integrable densities; when the domain $E$ is compact, the subscript $c$ will be ignored.

For $0<s<1$, the space of bounded $s$-Hölder continuous functions on $E$ with $s$-Hölder coefficient less than $H$ is denoted  by $C^{s}_H(E)$, i.e. ,
\[ C^s_H(E) := \left\{f: E \to \mathbb{C} :\sup_{\mathbf{x},\mathbf{y} \in E,\mathbf{x}\neq \mathbf{y}} \frac{| f(\mathbf{x}) - f(\mathbf{y}) |}{|\mathbf{x}-\mathbf{y}|^s}\leqslant H \right\}. \] 
The corresponding homogeneous space, i.e.,  the quotient space of $C^{s}_H(E)$ obtained by identifying functions differing by a constant, is denoted by $\dt{C}^{s}_H(E)$. Likewise, we define the $s$-Hölder space $C^s(E)\vcentcolon= \cup_{H>0}C_H^s(E)$ and use $\dt{C}^s(E)$ to denote its homogeneous space. In many situations, the notation for the homogeneous space is also used to denote a collection of representatives from each equivalence class satisfying $f(\mathbf{x_0}) = \mathbf{0}$ with a $\mathbf{x_0}\in E$ fixed; we will call it a realized homogeneous space in this case. Note that, under this interpretation, $\dt{C}^s(E)$ becomes a normed space with the original $s$-Hölder seminorm and its topological dual $(\dt{C}^s(E))^{\ast}$ is thus well-defined.

We write $\int_{E}\ \cdot \ d\mathbf{x}$ when integrating with respect to the $n$-dimensional Lebesgue measure on $E$. For $p>0$, the space of Lebesgue $p$-th integrable functions on $E$ is written as $L^p(E)$ and the space of $p$-th summable sequences on a set $X$ is written as $\ell^p(X)$, whose elements are expressed as $(a_i)_{i\in X}$. Moreover, when the underlying ground space is clear, they are sometimes abbreviated as $L^p$ and $\ell^p$.
\subsection{$s$-Wasserstein distance}\label{subsec2.2}
We begin by recalling several basic properties of the $s$-Wasserstein distance. For further background, see Santambrogio \cite{bib3}. For $E\subseteq \mathbb{R}^n$ and $E'\subseteq \mathbb{R}^n$ closed, given two probability measures $\mu\in\mathcal{P}(E)$, $\nu\in\mathcal{P}(E')$, a Borel measurable map $T$: $E \rightarrow E'$ is called a \textbf{transport map} from $\mu$ to $\nu$ if $\nu(A) = \mu(T^{-1}(A))$ for all Borel sets $A\subseteq E'$; in this case, we call $\nu$ the pushforward of $\mu$ under $T$ and denote it as $\nu = T_{\#}\mu$. The set of all \textbf{transport plans} from $\mu$ to $\nu$ is then defined by
$$\Gamma(\mu, \nu) \vcentcolon  = \{\gamma \in \mathcal{P}(E \times E') \ | \ \pi^1_{\#}\gamma = \mu,  \pi^2_{\#}\gamma = \nu \},$$
where $\pi^1 : E \times E' \to E$ and $\pi^2 : E \times E' \to E'$ are the canonical projections. We now define the general $s$-Wasserstein distance for $0<s\leqslant1$; in the case of $s=1$, it is informally also called the earth mover's distance (EMD).\\

\noindent \textbf{Definition 2.1}\phantomsection\label{Definition 2.1} ($s$-Wasserstein distance for $0<s\leqslant1$) For $\mu$, $\nu\in\mathcal{P}_s(E)$ and $0<s\leqslant1$, the \textbf{$\pmb{s}$-Wasserstein distance} $W_s(\mu, \nu)$ is defined as \begin{equation}
W_s(\mu, \nu) \vcentcolon  = \inf_{\gamma \in \Gamma(\mu, \nu)}\int_{E\times E}|\mathbf{x}-\mathbf{y}|^sd\gamma(\mathbf{x},\mathbf{y}).\tag{eq2.1}\label{eq2.1}\end{equation}

\noindent The minimizer is always attained by some $\gamma_{\text{min}}$ \cite[Theorem 1.7]{bib3}, and is called the \textbf{optimal transport plan} between $\mu$ and $\nu$.\\

\noindent \textsc{Remark 2.2} (Choice of the power $s$) In many standard texts, the $s$-Wasserstein distance is defined for $s\geqslant1$ with an extra power $1/s$ on the right-hand side of (\ref{eq2.1}), to ensure it as a distance. Besides the interest in the concave ground distances, we focus on the cases of $0<s\leqslant1$ due to the fact that the associated dual problem is well-characterized in terms of wavelets.\\

\noindent \textsc{Remark 2.3}\phantomsection\label{Remark 2.3} (Invariance of the distance) From the definition, it is apparent that $W_s(\mu,\nu)$ is invariant whether $\mu,\nu$ are considered to be in the space $\mathcal{P}_s(E)$ or $\mathcal{P}_s(F)$, as long as $E$ and $F$ contain the support of two measures. By the homogeneity of the ground distance, translating $\mu$ and $\nu$ by the same vector also does not change $W_s(\mu,\nu)$. Those points will be implicitly used throughout the paper.\\

\noindent \textbf{Proposition 2.4}\phantomsection\label{Proposition 2.4} (Kantorovich-Rubinstein duality of $W_s$ for $0<s\leqslant1$) Given $\mu, \nu\in\mathcal{P}_s(E)$ and $0<s\leqslant1$, then
\begin{equation}W_s(\mu, \nu) = \sup_{f\in C^{s}_{1}(E)}\int_{E}f(\mathbf{x}) d(\mu-\nu)(\mathbf{x}),\notag\end{equation}
and the maximum is attained by some  $f_{\text{max}}$.

\noindent \emph{Proof}. See \cite[Proposition 3.1]{bib3} for the duality; see \cite[Theorem 1.39, Proposition 3.1]{bib3} for the existence of the maximizer.\qed \\

\noindent Lastly, we record a convergence result. For lack of a reference, we include a proof.\\

\noindent \textbf{Proposition 2.5}\phantomsection\label{Proposition 2.5} (Convergence from $W_s$ to $W_1$) If $\mu, \nu\in\mathcal{P}_1(E)$, then $$\lim_{s\to1^-}W_s(\mu, \nu) = W_1(\mu, \nu).$$
\emph{Proof}. By Jensen's inequality and Definition \hyperref[Definition 2.1]{2.1}, we see that, for all $s \in (0,1)$, \begin{equation}W_s(\mu,\nu) \leqslant W_1(\mu,\nu)^s.\tag{eq2.2}\label{eq2.2}\end{equation} On the other hand, for any sequence $s_n \to 1^-$, let $(\gamma_n)_{n\in\mathbb{N}}$ be a sequence of optimal transport plans between $\mu$ and $\nu$, each with respect to $W_{s_n}$. After interpreting $\mu$ and $\nu$ to be measures on $\mathbb{R}^n$, as in Remark \hyperref[Remark 2.3]{2.3}, the proof of \cite[Theorem 1.7]{bib3} shows that, up to a subsequence, there exists $\gamma_* \in \Gamma(\mu,\nu)$  so that $\gamma_n \rightharpoonup \gamma_*$ weakly in the duality with bounded, continuous functions on $\mathbb{R}^n \times \mathbb{R}^n$. Then, by Ambrosio et al. \cite[Lemma 5.1.7]{bib30}, for all open $B_r \subseteq \mathbb{R}^n\times \mathbb{R}^n$ with $r>0$, we have 
\begin{equation}
       \liminf_{n \to +\infty} \int_{B_r} |x-y| d \gamma_n(x,y) \geq \int_{B_r} |x-y| d \gamma_*(x,y)\tag{eq2.3}\label{eq2.3}.
\end{equation}
Furthermore, for all $r>0$,
\begin{align*}
\int_{\mathbb{R}^n \times \mathbb{R}^n} |x-y|^{s_n} d \gamma_n(x,y) &\geq \int_{B_{r} } |x-y|^{s_n} d \gamma_n(x,y) \\
& = \int_{B_{r} } |x-y|^{s_n}-|x-y| d \gamma_n(x,y) + \int_{B_{r} } |x-y| d \gamma_n(x,y).
\end{align*}
Using that $|x-y|^{s_n}-|x-y| \to 0$ locally uniformly as $s_n\to1^-$ and (\ref{eq2.3}), we obtain
\[ \liminf_{n \to +\infty} \int_{\mathbb{R}^n \times \mathbb{R}^n} |x-y|^{s_n} d \gamma_n(x,y) \geq \int_{ B_r} |x-y| d \gamma_*(x,y).\]
Sending $r \to +\infty,$ the monotone convergence theorem ensures that 
\[ \liminf_{n \to +\infty} \int_{\mathbb{R}^n \times \mathbb{R}^n} |x-y|^{s_n} d \gamma_n(x,y) \geq \int_{\mathbb{R}^n \times \mathbb{R}^n} |x-y| d \gamma_*(x,y).\]
Since $\gamma_n$ were optimal transport plans, together with (\ref{eq2.2}), we have
\[ \int_{\mathbb{R}^n \times \mathbb{R}^n} |x-y| d \gamma_*(x,y) \leq \liminf_{n \to +\infty} W_{s_n}(\mu,\nu) \leq \liminf_{n \to +\infty}W_1(\mu,\nu)^{s_n} = W_1(\mu,\nu). \]
Due to the fact that $\gamma_* \in \Gamma(\mu,\nu)$, the left hand side must be larger than $W_1(\mu,\nu)$. This shows that equality holds throughout, so $\liminf_{n \to +\infty} W_{s_n}(\mu,\nu) = W_1(\mu,\nu)$. Finally, under (\ref{eq2.2}) again, we also have $\limsup_{n \to +\infty} W_{s_n}(\mu,\nu) \leqslant W_1(\mu,\nu)$. They imply $\lim_{s \to 1^-} W_s(\mu,\nu) = W_1(\mu,\nu)$.\qed

\subsection{Orthonormal wavelet systems}\label{subsec2.3}
We now recall background on orthonormal wavelet systems; for further materials on wavelets, see Meyer's book \cite{bib4}. The orthonormal wavelet system we need comes naturally from a \textbf{multiresolution approximation} on $L^2(\mathbb{R}^n)$, which is an increasing sequence of closed subspaces $\{V_{j}\}_{j\in\mathbb{Z}}$ in $L^2(\mathbb{R}^n)$ such that 
\begin{itemize}[leftmargin = 1.5em, label = {}]

\item (1) $$\bigcap_{j\in\mathbb{Z}} V_j = \{0\} \qquad  \text{and} \qquad \overline{\bigcup_{j\in\mathbb{Z}} V_j} = L^2(\mathbb{R}^n),$$\noindent where the bar indicates the closure.

\item (2) For all $f \in L^2(\mathbb{R}^n)$, $j\in\mathbb{Z}$, and $\mathbf{k}\in\mathbb{Z}^n$,
$$f(\mathbf{x}) \in V_j \Leftrightarrow f(2\mathbf{x}) \in V_{j+1} \qquad
\text{and} \qquad
f(\mathbf{x}) \in V_{0} \Leftrightarrow f(\mathbf{x} - \mathbf{k}) \in V_{0}$$

\item (3) There exists a \textbf{scaling function} $\phi \in V_0$ such that it integrates to 1 and $\{\phi_{(j,\mathbf{k})}\}_{\mathbf{k}\in\mathbb{Z}^n}$ forms an orthonormal
Hilbert basis of $V_j$ for all $j\in\mathbb{Z}$, where
$$\phi_{(j,\mathbf{k})}(\mathbf{x}) \vcentcolon  = 2^{\frac{nj}{2}}\phi(2^j\mathbf{x} - \mathbf{k}).$$
\end{itemize}
With a multiresolution approximation $\{V_{j}\}_{j\in\mathbb{Z}}$, we can denote $W_j$ to be the orthogonal
complement of $V_j$ in $V_{j+1}$. Then, by definition, 
\begin{equation}
L^2(\mathbb{R}^n) = V_{j_0}\bigoplus \left(\bigoplus_{j\geqslant j_0} W_j\right) = \bigoplus_{j\in\mathbb{Z}} W_j\tag{eq2.4}\label{eq2.4}
\end{equation}
for all $j_0\in\mathbb{Z}$, where $\bigoplus$ is the Hilbert direct sum. Note that the summands in (\ref{eq2.4}), whether in the middle sum or the right sum, are all mutually orthogonal. The definition of an orthonormal wavelet system is as follows.\\

\noindent \textbf{Definition 2.6}\phantomsection\label{Definition 2.6}
(Orthonormal wavelet system) Given the notations above, the collection of subspaces $\{V_{j}\}_{j\in\mathbb{Z}},\{W_{j}\}_{j\in\mathbb{Z}}$ of $L^2(\mathbb{R}^n)$ is called an \textbf{orthonormal wavelet system} if there exists a set of functions $\{\psi_{\epsilon}\}_{\epsilon\in\alpha}\in W_0$ indexed by a finite set $\alpha$ such that each of them integrates to 0 and $\{\psi_{(\epsilon,(j,\mathbf{k}))}\}_{\epsilon\in\alpha,\mathbf{k}\in\mathbf{Z}^n}$ is an orthonormal Hilbert basis of $W_j$ for all $j\in\mathbb{Z}$, where
$$\psi_{(\epsilon,(j,\mathbf{k}))}(\mathbf{x})\vcentcolon  =2^{\frac{nj}{2}}\psi_{\epsilon}(2^j\mathbf{x}-\mathbf{k}).$$
The functions $\{\psi_{\epsilon}\}_{\epsilon\in\alpha}$ are called the \textbf{wavelet functions}.\\

Some authors require more conditions on $\phi$ and $\{\psi_{\epsilon}\}_{\epsilon\in\alpha}$ for them to be qualified as parts of wavelet systems. The ones
in Definition \hyperref[Definition 2.6]{2.6} are sufficient for us.

It is clear that an orthonormal wavelet system is completely determined by its scaling function $\phi$ and its wavelet functions $\{\psi_{\epsilon}\}_{\epsilon\in\alpha}$. So, under a slight abuse of notation, we shall denote a system by $(\phi,\psi)$. The existence of orthonormal wavelet systems is not obvious, and a detailed construction can be found in \cite[Chapter 3]{bib4}. In fact, as shown in Meyer's book, it suffices to find a system $(\phi,\psi)$ for $n=1$ dimension. In this case, from a constructive perspective, the index set $\alpha$ is always a singleton, and we will assume all wavelets on $\mathbb{R}$ are like this. Then, using what is known as the tensor product method, one can build a system for general dimension $n\in\mathbb{N}$. In $n=2$, for example, let $$\phi_{(j,\mathbf{k})}^{(2)}(\mathbf{x}) \vcentcolon  = 2^{\frac{j}{2}}\phi(2^jx_1-k_1)2^{\frac{j}{2}}\phi(2^jx_2-k_2) \qquad (j, \mathbf{k})\in\mathbb{Z}\times\mathbb{Z}^2$$
and
$$\begin{aligned}
\psi_{(1, (j,\mathbf{k}))}^{(2)}(\mathbf{x}) &\vcentcolon  = 2^{\frac{j}{2}}\phi(2^jx_1-k_1)2^{\frac{j}{2}}\psi(2^jx_2-k_2),\\
\psi_{(2, (j,\mathbf{k}))}^{(2)}(\mathbf{x}) &\vcentcolon  = 2^{\frac{j}{2}}\psi(2^jx_1-k_1)2^{\frac{j}{2}}\psi(2^jx_2-k_2), \qquad (j, \mathbf{k})\in\mathbb{Z}\times\mathbb{Z}^2.\\
\psi_{(3, (j, \mathbf{k}))}^{(2)}(\mathbf{x}) &\vcentcolon  = 2^{\frac{j}{2}}\psi(2^jx_1-k_1)2^{\frac{j}{2}}\phi(2^jx_2-k_2)
\end{aligned}$$
Then $\phi_{(0,\mathbf{0})}^{(2)}$ is the initial scaling function and $\{\psi_{(\epsilon,(0,\mathbf{0}))}^{(2)}\}_{\epsilon\in\llbracket 3 \rrbracket}$ are the initial wavelet functions for the new system. For more general $n\in \mathbb{N}$, proceed inductively and obtain the collections $\{\phi_{(j,\mathbf{k})}^{(n)}\}_{j\in\mathbb{Z}, \mathbf{k}\in\mathbb{Z}^n}$ and $\{\psi_{(\epsilon,(j,\mathbf{k}))}^{(n)}\}_{\epsilon\in\llbracket 2^n-1 \rrbracket, j\in\mathbb{Z},\mathbf{k}\in\mathbb{Z}^n}$; note that due to combinatorics, $\epsilon$ has to range from 1 to $2^{n}-1$. When $n$ is fixed, we will ignore the superscripts in the following sections, even if it is constructed from $n=1$ like above.

Now, given an orthonormal wavelet system $(\phi,\psi)$ for $L^2(\mathbb{R}^n)$, the relationships (\ref{eq2.4}) explicitly mean: for any $f\in L^2(\mathbb{R}^n)$ and $j_0\in\mathbb{Z}$, 
\begin{equation}
f(\mathbf{x}) = \sum_{\mathbf{k}\in\mathbb{Z}^n}f_{(j_0, \mathbf{k})}^{(a)}\phi_{(j_0, \mathbf{k})}(\mathbf{x}) + \sum_{\substack{\epsilon\in\alpha,\\ j\geqslant j_0,\  \mathbf{k}\in\mathbb{Z}^n}}f_{(\epsilon, (j, \mathbf{k}))}^{(d)}\psi_{(\epsilon, (j, \mathbf{k}))}(\mathbf{x}),\tag{eq2.5}\label{eq2.5}
\end{equation}
for some complex-valued coefficients $f_{(j_0,\mathbf{k})}^{(a)}$ and $f_{(\epsilon,(j,\mathbf{k}))}^{(d)}$, and (\ref{eq2.5}) is called the \textbf{scaling expansion} of $f$ \textbf{at level $\pmb{j_0}$}; alternatively, we have 
\begin{equation}
f(\mathbf{x}) =  \sum_{\substack{\epsilon\in\alpha,\\ j\in\mathbb{Z},\  \mathbf{k}\in\mathbb{Z}^n}}f_{(\epsilon, (j, \mathbf{k}))}^{(d)}\psi_{(\epsilon, (j, \mathbf{k}))}(\mathbf{x}).\tag{eq2.6}\label{eq2.6}
\end{equation}
for some complex-valued coefficients $f_{(\epsilon,(j,\mathbf{k}))}^{(d)}$, and (\ref{eq2.6}) is called the \textbf{full wavelet expansion} of $f$. The equalities are typically understood as convergence in  $L^2(\mathbb{R}^n)$.

The coefficient $f_{(j,\mathbf{k})}^{(a)}$ is called the \textbf{approximation coefficient} of $f$ \textbf{at level $\pmb{j}$ and translation $\mathbf{k}$}, while the $f_{(\epsilon,(j,\mathbf{k}))}^{(d)}$ is called the \textbf{detail coefficient} of $f$ \textbf{at level $\pmb{j}$ and translation $\mathbf{k}$}. Sometimes they are collectively called the wavelet coefficients, under abuse of terminology. Due to the orthonormality of our wavelet basis, they can be computed as
\begin{align}
f_{(j,\mathbf{k})}^{(a)} &= \int_{\mathbb{R}^n}f(\mathbf{x})\phi_{(j,\mathbf{k})}^{\ast}(\mathbf{x})d\mathbf{x}\tag{eq2.7}\label{eq2.7},\\
f_{(\epsilon, (j,\mathbf{k}))}^{(d)} &= \int_{\mathbb{R}^n}f(\mathbf{x})\psi_{(\epsilon, (j,\mathbf{k}))}^{\ast}(\mathbf{x})d\mathbf{x}\tag{eq2.8}\label{eq2.8},
\end{align}
where the $\ast$ denotes the complex conjugate. More generally, for $f\notin L^2(\mathbb{R}^n)$, we will still use $f_{(j,\mathbf{k})}^{(a)}$ and $f_{(\epsilon,(j,\mathbf{k}))}^{(d)}$ to represent the values of (\ref{eq2.7}) and (\ref{eq2.8}) as long as the integrals are well-defined. Of course, the convergence of the formal expansions like (\ref{eq2.5}) and (\ref{eq2.6}) need to be analyzed differently.

Lastly, we call an orthonormal wavelet system $(\phi,\psi)$ to be of \textbf{regularity $\pmb{r}$} if $\phi$ and $\{\psi_{\epsilon}\}_{\epsilon\in\alpha}$ are continuously differentiable up to order $r$ and all their derivatives are decaying faster than the reciprocal of any polynomial at infinity. We call the wavelet system \textbf{compactly supported} if  $\phi$ and $\{\psi_{\epsilon}\}_{\epsilon\in\alpha}$ are. See \cite[Chapter 3]{bib4} for the existence of $r$ regular, compactly supported orthonormal wavelet systems for any $r\in\mathbb{N}$.

\subsection{Discrete wavelet transform}\label{subsec2.4}
In view of (\ref{eq2.7}) and (\ref{eq2.8}), to get an explicit expansion of a function in $L^2(\mathbb{R}^n)$, as in (\ref{eq2.5}) or (\ref{eq2.6}), it amounts to computing many integrals. Alternatively, given an initialization of finitely many, consecutive approximation coefficients at a level, certain coefficients below the initialization level can be computed in linear time with a cascaded pyramidal algorithm called \textbf{the discrete wavelet transform} (\textbf{DWT}); see Mallat \cite[Chapter 7]{bib6}. This is essentially due to the structure of the multiresolution approximation.

We focus on the DWT for $n=1$.\footnote{It is worth mentioning that, for our purpose, the discrete wavelet transform is used to compute the actual approximation/detail coefficients of some \textit{non-discrete functions}, so we frame the transform in a way to reflect this connection. On the other hand, it is certainly possible to frame (and treat) the discrete wavelet transform as an algorithm fully in a discrete universe, independent of any connection to its continuous counterpart. The more well-known discrete Fourier transform also has such two perspectives.} Given an orthonormal wavelet system $(\phi,\psi)$ and $f\in L^2(\mathbb{R})$, the algorithm takes an input approximation vector 
\begin{equation}
\mathbf{A}^{j_0+M} = \left(f_{(j_0+M,0)}^{(a)}, f_{(j_0+M,1)}^{(a)}, \cdots ,f_{(j_0+M,2^{M}-1)}^{(a)}\right)\tag{eq2.9}\label{eq2.9}
\end{equation}
of length $2^M$ for some $M\in\mathbb{N}$ and $j_0\in\mathbb{Z}$, consisting of a consecutive collection of approximation coefficients of $f$ at a level $j_0+M$; then, it gives
\begin{align}
f_{(j-1,k)}^{(a)} &= A^{j-1}_{k} \vcentcolon  = \sum_{l=N_{l,\phi}}^{N_{u,\phi}}A^{j}_{l}g_{l-2k},\tag{eq2.10}\label{eq2.10}\\
f_{(j-1,k)}^{(d)} &= D^{j-1}_{k} \vcentcolon  = \sum_{l=N_{l,\psi}}^{N_{u,\psi}}A^{j}_{l}h_{l-2k}, \tag{eq2.11}\label{eq2.11}
\end{align}
with
$$\quad j = j_0+M,\cdots,j_0+1, \quad k = 0,\cdots, 2^{j-j_0-1}-1,$$
where the vectors $\mathbf{g}$ and $\mathbf{h}$ are the filters associated with $\phi_{(0,0)}$ and $\psi_{(0,0)}$ respectively, and their lengths (bounded by $N_{l,\phi}$, $N_{u,\phi}$, $N_{l,\psi}$, $N_{u,\psi}$) are fully determined by $\phi_{(0,0)}$ and $\psi_{(0,0)}$. Iteratively, the algorithm computes an approximation vector $\mathbf{A}^{j-1}$ and a detail vector $\mathbf{D}^{j-1}$ with halved length at the level $j-1$ given the approximation vector $\mathbf{A}^j$ at level $j$. At last, the discrete wavelet transform outputs a vector $\mathbf{O}^{j_0,M}$ of length $2^M$, with the first element being the single approximation coefficient $f_{(j_0,0)}^{(a)}$ at level $j_0$, then followed by the single detail coefficient $f_{(j_0,0)}^{(d)}$ at level $j_0$, then the double detail coefficients $f_{(j_0+1,0)}^{(d)}$ and $f_{(j_0+1,1)}^{(d)}$ at level $j_0+1$, up to the detail coefficients of length $2^{M-1} $ at level $j_0+M-1$. More explicitly,
\begin{equation}
\mathbf{O}^{j_0,M} = \bigl(f_{(j_0,0)}^{(a)},f_{(j_0,0)}^{(d)}, f_{(j_0+1,0)}^{(d)}, f_{(j_0+1,1)}^{(d)},\cdots,f_{(j_0+M-1,0)}^{(d)},\cdots,f_{(j_0+M-1,2^{M-1}-1)}^{(d)}\bigr).\tag{eq2.12}\label{eq2.12}
\end{equation}
\textsc{Remark 2.7}\phantomsection\label{Remark 2.7} (Modes for DWT) The algorithm introduced should be technically called the DWT  with the periodic extension mode, meaning that the approximation vector at each iteration is automatically extended in a periodic fashion on the boundary. This extension is necessary to make the sums in (\ref{eq2.10}) and (\ref{eq2.11}) well-defined when the coefficients to be computed are near two edges. There are other extension modes, such as the zeros extension mode we will use in Section \ref{sec4}, which simply extends the approximation vector at each level by adding zeros on both sides. The reason we choose to start with the periodic extension mode is that the output vector $\mathbf{O}^{j_0,M}$ has the same length as the input vector $\mathbf{A}^{j_0,M}$ and is more suitable to understand core structure of the cascaded pyramidal algorithm. In other cases, depending on the length of the filters, a few more coefficients will be computed in each iteration level, resulting in the final output vector $\mathbf{O}^{j_0,M}$ being slightly longer than $2^M$; the explicit output is slightly more complicated to write due to indices, but the logic is the same.\\

\noindent \textsc{Remark 2.8} (Generality of the input vector) The most general input approximation vector should be in the form of 
\begin{equation}
\mathbf{A}^{j_0+M,k_0} = \left(f_{(j_0+M,k_0)}^{(a)}, f_{(j_0+M,1+k_0)}^{(a)}, \cdots ,f_{(j_0+M,2^{M}+k_0-1)}^{(a)}\right)\tag{eq2.12}\label{eq2.12}
\end{equation}
of length $2^M$ for some $M\in\mathbb{N}$, $j_0\in\mathbb{Z}$, and $k_0\in\mathbb{Z}$. However, one can always shift the function $f$ to some $\tilde{f}$ so that, by using a change of variables in (\ref{eq2.7}), (\ref{eq2.12}) equals to 
 $$\left(\tilde{f}_{(j_0+M,0)}^{(a)}, \tilde{f}_{(j_0+M,1)}^{(a)}, \cdots ,\tilde{f}_{(j_0+M,2^{M}-1)}^{(a)}\right).$$
 Therefore the one in (\ref{eq2.9}) can be viewed as general, and the reason for doing so is to make the formulas above and below cleaner.
 
 In fact, as one shall see in Section \ref{sec4} when we actually implement DWT, the specific scenario allows us to directly use (\ref{eq2.9}) without even considering the specific $\tilde{f}$ being shifted to. \\

\noindent Now, to obtain certain approximation and detail coefficients below a level, one can obtain a proper input approximation vector on that level and use the DWT. However, computationally, the theoretical way of purely using (\ref{eq2.7}) and (\ref{eq2.8}) and the one with the DWT amount to the same number of integrals to be computed due to the initialization required in the later one. The practical usefulness of the DWT lies in making some approximations to simplify the initialization. 

One of the simplest approximation goes as follows: assume the system $(\phi,\psi)$ initially chosen is compactly supported around the origin (or decay quickly), then for $j_0+M$ positive and large and $f$ smooth enough, $f$ is almost constant on the domain where $\phi_{(j_0+M,k)}(\mathbf{x})$ is nonzero, with support roughly at $2^{-(j_0+M)}k$. Then by (\ref{eq2.7}) and the fact that $\phi$ integrates to 1,
\begin{equation}f_{(j_0+M,k)}^{(a)}\approx 2^{-\frac{j_0+M}{2}}f(2^{-(j_0+M)}k).\tag{eq2.13}\label{eq2.13}
\end{equation}
There exist other more sophisticated approximation schemes, particularly the one in Zhang et al. \cite{bib8}, where the (\ref{eq2.13}) we mentioned corresponds to the most major term in their method. We experimented and found that those two approximations of initialization produce almost no difference when $j_0+M$ is large, which is the natural choice we shall take.

The whole development, from DWT to various approximations, can be carried over to the case for $n>1$. This is in principle due to that the algorithm by Mallat \cite{bib6} is based on the multiresolution approximation, which makes no difference in higher dimensions. The same goes for the approximations shown above; we will not detail here, as in the numerics section (Section \ref{sec4}), we focus on the case for $n=1$ to illustrate our result more transparently.

\section{Main results}\label{sec3}
We now turn to our main results. First, we introduce the original wavelet $s$-Wasserstein distance in \cite{bib1} and explain the problems with this original formulation. Next, we consider two approaches to fixing the problem: a naive solution, which is close in spirit to the original formulation, and a novel reformulation, from a slightly different starting point. We discuss the benefits of our novel reformulation and conclude by describing how it provides a linear embedding, with an equivalent distance to $W_s$ (up to an error).

We begin by recalling some results about wavelet characterizations of H\"older spaces, which are the foundation of the various wavelet $s$-Wasserstein distances we consider.
\subsection{Wavelets, $s$-Hölder spaces, and the $s$-Wasserstein distance}\label{subsec3.1}
Throughout the following presentation, the dimension $n$ and the distance parameter $s\in(0,1)$ are fixed, unless otherwise specified. We assume the orthonormal wavelet system $(\phi,\psi)$ to be 1-regular and compactly supported.

Let $d\mu = p_{\mu}(\mathbf{x})d\mathbf{x}$, $d\nu = p_{\nu}(\mathbf{x})d\mathbf{x}\in\mathcal{P}_{c,\text{ab2}}(\mathbb{R}^n)\subseteq\mathcal{P}_{s}(\mathbb{R}^n)$. Based on Proposition \hyperref[Proposition 2.4]{2.4}, the $s$-Wasserstein distance between them is 
\begin{equation}
W_s(\mu, \nu) = \sup_{f\in C^{s}_{1}(\mathbb{R}^n)\cap L^2(\mathbb{R}^n)}\int_{\mathbb{R}^n}f(\mathbf{x}) (p_{\mu}(\mathbf{x})-p_{\nu}(\mathbf{x}))d\mathbf{x},\notag\end{equation}
where we can impose the extra constraint of $f\in L^2(\mathbb{R}^n)$ since $\mu$ and $\nu$ are compactly supported. Furthermore, since $f$, $p_{\mu}-p_{\nu}\in L^2(\mathbb{R}^n)$, by Parseval's theorem and (\ref{eq2.5}), for any $j_0\in\mathbb{Z}$, we have
\begin{equation}
W_s(\mu,\nu) = \sup_{\substack{f \in C^{s}_{1}(\mathbb{R}^n)\\ \cap L^2(\mathbb{R}^n)}} \Biggl( \sum_{\mathbf{k} \in \mathbb{Z}^n} f_{(j_0,\mathbf{k})}^{(a)} (p_{\mu}-p_{\nu})_{(j_0,\mathbf{k})}^{(a)} +  \notag
\sum_{\substack{\epsilon \in \alpha,\\ j \geqslant j_0, \mathbf{k} \in \mathbb{Z}^n}} f_{(\epsilon, (j, \mathbf{k}))}^{(d)}(p_{\mu}-p_{\nu})_{(\epsilon, (j, \mathbf{k}))}^{(d)} \Biggr),\notag
\end{equation}
and likewise, by (\ref{eq2.6}), we have
\begin{equation}
W_s(\mu,\nu) = \sup_{f\in C^{s}_{1}(\mathbb{R}^n)\cap L^2(\mathbb{R}^n)}\sum_{\substack{\epsilon\in\alpha,\\ j\in\mathbb{Z},\mathbf{k}\in\mathbb{Z}^n}}f_{(\epsilon, (j,\mathbf{k}))}^{(d)}(p_{\mu}-p_{\nu})_{(\epsilon, (j,\mathbf{k}))}^{(d)}.\notag
\end{equation}
As we will see, it will be convenient to allow the $s$-Hölder constant $H>0$ in the constraint to vary. Note   that the $s$-Wasserstein distance scales with $H$ as follows:
\begin{equation}
H W_s(\mu, \nu) = \sup_{f\in C^{s}_{H}(\mathbb{R}^n)\cap L^2(\mathbb{R}^n)}\int_{\mathbb{R}^n}f(\mathbf{x}) (p_{\mu}(\mathbf{x})-p_{\nu}(\mathbf{x}))d\mathbf{x}.\tag{eq3.1}\label{eq3.1}\end{equation}
Therefore, for any $H>0$ and $j_0\in\mathbb{Z}$, computing $H W_s(\mu,\nu)$ reduces to an optimization problem over the set
\begin{align}\text{Wav}_{j_0}^{a}(C^{s}_{H}(\mathbb{R}^n)) \vcentcolon  = \Bigl{\{}&(f_{(j_0,\mathbf{k})}^{(a)})_{\mathbf{k}\in\mathbb{Z}^n}\bigoplus (f_{(\epsilon, (j,\mathbf{k}))}^{(d)})_{\substack{\epsilon \in \alpha, \\j \geqslant j_0, \mathbf{k} \in \mathbb{Z}^n}}\in \notag \ell^2\ \Big|\ f \in C^{s}_{H}(\mathbb{R}^n) \Bigr{\}},\notag
\end{align}
or the set
\begin{equation}
\text{Wav}^{d}(C^{s}_{H}(\mathbb{R}^n)) \vcentcolon  = \left\{ (f_{(\epsilon, (j,\mathbf{k}))}^{(d)})_{\substack{\epsilon\in\alpha,\\j\in\mathbb{Z},\mathbf{k}\in\mathbb{Z}^n}}\in \ell^2\ \Big|\ f \in C^{s}_{H}(\mathbb{R}^n) \right\}.\notag
\end{equation}

In \cite{bib1}, the authors choose to work with the set $\text{Wav}_{0}^{a}(C^s_{H}(\mathbb{R}^n))$, derived from the scaling expansion at level 0, due to the following theorem.\\

\noindent \textbf{Theorem 3.1}\phantomsection \label{Theorem 3.1} (Characterizations of $C^s(\mathbb{R}^n)$ via wavelet systems \cite[Theorem 6.4.5]{bib4}) A function $f\in L_{\text{loc}}^1(\mathbb{R}^n)$ belongs to $C^s(\mathbb{R}^n)$ if and only if, when in the scaling expansion at level 0, its wavelet coefficients satisfy
\begin{equation}
|f_{(0,\mathbf{k})}^{(a)}|\leqslant C_0 \qquad \mathbf{k}\in\mathbb{Z}^n,
\notag
\end{equation}
for some constant $C_0\geqslant 0$ and
\begin{equation}
|f_{(\epsilon,(j,\mathbf{k}))}^{(d)}|\leqslant C_12^{-j(s+\frac{n}{2})} \qquad \epsilon\in\alpha,j\geqslant0,\mathbf{k}\in\mathbb{Z}^n,
\notag
\end{equation}
for some constant $C_1\geqslant0$.\\

\noindent This theorem shows that the $s$-Hölder regularity of a function is characterized by the rate of decay of its wavelet coefficients. However, it is qualitative, in the sense that it does not explicitly give a connection between the $s$-Hölder constant $H$ of a function and the constants $C_0$ and $C_1$ in the bounds of its wavelet coefficients. One of the main contributions of \cite{bib1} is the following quantitative version of the above theorem, which is the foundation of the wavelet $s$-Wasserstein distances. We state this result in an expanded way for the reader's convenience.\\
 
\noindent \textbf{Lemma 3.2}\phantomsection\label{Lemma 3.2} (Quantitative characterizations of $C^s(\mathbb{R}^n)$ via wavelet systems \cite[Lemma 1 \& Appendix]{bib1}) Fix $H>0$. For all $f\in C_{H}^s(\mathbb{R}^n)$, when in the scaling expansion at level 0, the wavelet coefficients are bounded by 
\begin{equation}
|f_{(0,\mathbf{k})}^{(a)}|\leqslant ||f||_{L^{\infty}(\mathbb{R}^n)}+\frac{H}{a_{11}(\phi,s)}, \qquad \mathbf{k}\in\mathbb{Z}^n,
\notag
\end{equation}
and
\begin{equation}
|f_{(\epsilon,(j,\mathbf{k}))}^{(d)}|\leqslant \frac{H2^{-j(s+\frac{n}{2})}}{a_{12}(\psi,s)}, \qquad \epsilon\in\alpha,j\geqslant0,\mathbf{k}\in\mathbb{Z}^n,
\notag
\end{equation}
where $a_{11}(\phi,s), a_{12}(\psi,s)$ are constants depending on the $s$ and the wavelet system $(\phi,\psi)$, via
\begin{align}
a_{11}(\phi,s)&\vcentcolon  =\frac{1}{\inf_{\textbf{r}\in\mathbb{R}^n}\int_{\mathbb{R}^n}|\mathbf{x}-\mathbf{r}|^s|\phi(\mathbf{x})|d\mathbf{x}},\notag\\
a_{12}(\psi,s)&\vcentcolon  =\frac{1}{\max_{\epsilon\in\alpha}\inf_{\textbf{r}\in\mathbb{R}^n}\int_{\mathbb{R}^n}|\mathbf{x}-\mathbf{r}|^s|\psi_{\epsilon}(\mathbf{x})|d\mathbf{x}}.\notag
\end{align}

Conversely, given a function $\tilde{f}\in L_{\text{loc}}^1(\mathbb{R}^n)$ such that the wavelet coefficients are bounded by
\begin{equation}
|\tilde{f}_{(0,\mathbf{k})}^{(a)}|\leqslant C_0, \qquad \mathbf{k}\in\mathbb{Z}^n,\notag
\end{equation}
for some constant $C_0\geqslant 0$ and
\begin{equation}
|\tilde{f}_{(\epsilon,(j,\mathbf{k}))}^{(d)}|\leqslant C_12^{-j(s+\frac{n}{2})}, \qquad \epsilon\in\alpha,j\geqslant0,\mathbf{k}\in\mathbb{Z}^n,
\notag
\end{equation}
for some $C_1\geqslant0$, then 
\begin{equation}
\sum_{\mathbf{k}\in\mathbb{Z}^n}\tilde{f}_{(0, \mathbf{k})}^{(a)}\phi_{(0, \mathbf{k})}(\mathbf{x})+ \sum_{\substack{\epsilon\in\alpha,\\ j\geqslant0,\  \mathbf{k}\in\mathbb{Z}^n}}\tilde{f}_{(\epsilon, (j, \mathbf{k}))}^{(d)}\psi_{(\epsilon, (j, \mathbf{k}))}(\mathbf{x})\notag
\end{equation}
converges uniformly to $\tilde{f}$ on compact subsets. Furthermore, we have $\tilde{f}\in C_{\tilde{H}}^s(\mathbb{R}^n)$ for   $\tilde{H} \vcentcolon=  a_{21}(\phi,s)C_0+a_{22}(\psi,s)C_1$, where, following the notation in \cite{bib1},
\begin{equation}
a_{21}(\phi,s)\vcentcolon  = ||\textstyle\sum\phi^s||_{\infty}\vcentcolon  = \sup_{\substack{\mathbf{x}\neq\mathbf{y}\\\eta_k\in\{+1,-1\}}}\left|\frac{\sum_{\mathbf{k}\in\mathbb{Z}^n}\eta_k\phi_{(0,\mathbf{k})}(\mathbf{x})-\sum_{\mathbf{k}\in\mathbb{Z}^n}\eta_k\phi_{(0,\mathbf{k})}(\mathbf{y})}{|\mathbf{x}-\mathbf{y}|^s} \right|,\notag
\end{equation}
\\
\begin{equation}
||\textstyle\sum\psi||_{\infty}\vcentcolon  = \sup_{\substack{\mathbf{x}\in\mathbb{R}^n\\\eta_k\in\{+1,-1\}}}\left|\sum_{\mathbf{k}\in\mathbb{Z}^n,\epsilon\in\alpha}\eta_k\psi_{(\epsilon,(0,\mathbf{k}))}(\mathbf{x}) \right|,\notag
\end{equation}
\\
\begin{equation}
||\textstyle \sum\psi^{(i)}||_{\infty} \vcentcolon  = \sup_{\substack{\mathbf{x}\in\mathbb{R}^n\\\eta_k\in\{+1,-1\}}}\left|\sum_{\mathbf{k}\in\mathbb{Z}^n,\epsilon\in\alpha}\eta_k\frac{\partial\psi_{(\epsilon,(0,\mathbf{k}))}(\mathbf{x})}{\partial x_i} \right|,\notag
\end{equation}
\\
\begin{equation}
a_{22}(\psi,s) \vcentcolon  = \frac{\sum_{i=1}^n||\textstyle \sum\psi^{(i)}||_{\infty}}{2^{1-s}(2^{1-s}-1)}+\frac{2||\textstyle\sum\psi||_{\infty}}{1-2^{-s}}.\notag
\end{equation}\\

  The finiteness of all the constants appearing in the previous lemma follows easily from the regularity and compact support of the wavelet system.\\

\noindent \textsc{Remark 3.3}\phantomsection\label{Remark 3.3} (Characterizations of Lipschitz functions via wavelet systems) Theorem \hyperref[Theorem 3.1]{3.1} does not extend to the characterization of Lipschitz continuous functions in the case $s=1$. Instead, the functions with detail coefficients bounded by $2^{-j(1+n/2)}$ are those that belong to the Zygmund class, which is strictly larger than the space of Lipschitz continuous functions. The fact that the result does not extend to the case $s=1$ is also reflected by the fact that the constant  $a_{22}(\psi,s)$ cannot be continuously extended to $s=1$ in Lemma \hyperref[Lemma 3.2]{3.2}.\\

Next, we consider alternative bounds of the wavelet coefficients for functions in $C_H^s(\mathbb{R}^n)$. These do not appear in \cite{bib1} and are not part of the original formulation of the wavelet $s$-Wasserstein distance, but they will be useful later in our revised distances.\\

\noindent \textbf{Lemma 3.4} (Alternative bounds for one direction in Lemma \hyperref[Lemma 3.2]{3.2})\phantomsection\label{Lemma 3.4} Fix $H>0$. For all $f\in C_{H}^s(\mathbb{R}^n)$, when in the scaling expansion at level 0, the wavelet coefficients are bounded by
\begin{equation}
|f_{(0,\mathbf{k})}^{(a)}|\leqslant \frac{||f||_{L^{\infty}(\mathbb{R}^n)}}{a_{13}(\phi)}, \qquad \mathbf{k}\in\mathbb{Z}^n,
\tag{eq3.2}\label{eq3.2}
\end{equation}
and
\begin{equation}
|f_{(\epsilon,(j,\mathbf{k}))}^{(d)}|\leqslant \frac{H2^{-j(s+\frac{n}{2})}}{a_{14}(\psi,s)}, \qquad \epsilon\in\alpha,j\geqslant0,\mathbf{k}\in\mathbb{Z}^n,
\tag{eq3.3}\label{eq3.3}
\end{equation}
where $a_{13}(\phi), a_{14}(\psi,s)$ are defined as
\begin{align}
a_{13}(\phi)&\vcentcolon  =\frac{1}{||\phi||_{L^1(\mathbb{R}^n)}},\notag \qquad
a_{14}(\psi,s)\vcentcolon  =\frac{1}{\max_{\epsilon\in\alpha}||\psi_{\epsilon}||_{(\dt{C}^s(\mathbb{R}^n))^{\ast}}},\notag
\end{align}
where $||\cdot||_{(\dt{C}^s(\mathbb{R}^n))^{\ast}}$ is the operator norm on the topological dual space of $\dt{C}^{s}(\mathbb{R}^n)$.\\
\emph{Proof}. The first bound (\ref{eq3.2}) follows by definition (\ref{eq2.7}) and a change of variables:
$$|f_{(0,\mathbf{k})}^{(a)}|\leqslant\int_{\mathbb{R}^n}|f(\mathbf{x})\phi^{\ast}(\mathbf{x}-\mathbf{k})|d\mathbf{x}\leqslant||f||_{L^{\infty}(\mathbb{R}^n)}||\phi^{\ast}||_{L^1(\mathbb{R}^n)} = ||f||_{L^{\infty}(\mathbb{R}^n)}||\phi||_{L^1(\mathbb{R}^n)}.$$
The second bound (\ref{eq3.3}) follows by definition (\ref{eq2.8}) and a similar change of variables:
\begin{align}|f_{(\epsilon,(j,\mathbf{k}))}^{(d)}|&\leqslant\max_{\epsilon\in\alpha}\sup_{f\in C_{H}^{s}(\mathbb{R}^n)}\Biggl|\int_{\mathbb{R}^n}f(\mathbf{x})\psi_{(\epsilon, (j,\mathbf{k}))}^{\ast}(\mathbf{x})d\mathbf{x}\Biggr|\notag\\
&=\max_{\epsilon\in\alpha}\sup_{f\in C_{H}^{s}(\mathbb{R}^n)}\Biggl|\int_{\mathbb{R}^n}f^{\ast}(\mathbf{x})\psi_{(\epsilon, (j,\mathbf{k}))}(\mathbf{x})d\mathbf{x}\Biggr|\notag\\
&=\max_{\epsilon\in\alpha}\sup_{f\in C_{H}^{s}(\mathbb{R}^n)}\Biggl|\int_{\mathbb{R}^n}2^{-\frac{nj}{2}}f(2^{-j}(\mathbf{x}+\mathbf{k}))\psi_{\epsilon}(\mathbf{x})d\mathbf{x}\Biggr|\notag\\
&=2^{-\frac{nj}{2}}\max_{\epsilon\in\alpha}\sup_{f\in\dt{C}_{2^{-js}H}^{s}(\mathbb{R}^n)}\Biggl|\int_{\mathbb{R}^n}f(\mathbf{x}+\mathbf{k})\psi_{\epsilon}(\mathbf{x})d\mathbf{x}\Biggr|\notag\\
&=H2^{-j(s+\frac{n}{2})}\max_{\epsilon\in\alpha}\sup_{f\in\dt{C}_{1}^{s}(\mathbb{R}^n)}\Biggl|\int_{\mathbb{R}^n}f(\mathbf{x}+\mathbf{k})\psi_{\epsilon}(\mathbf{x})d\mathbf{x}\Biggr|\notag\\
&=H2^{-j(s+\frac{n}{2})}\max_{\epsilon\in\alpha}||\psi_{\epsilon}||_{(\dt{C}^s(\mathbb{R}^n))^{\ast}}\notag.
\end{align}
In the fourth equality, we replace $C_H^s(\mathbb{R}^n)$ with a realized $\dt{C}_H^s(\mathbb{R}^n)$, which is valid because all $\psi_{\epsilon}$ integrate to 0.\qed\\

\noindent \textsc{Remark 3.5} (Clarification for alternative bounds) 
We bound $|f_{(0,\mathbf{k})}^{(a)}|$ differently compared to Lemma \hyperref[Lemma 3.2]{3.2} since it leads to a simpler dependence on constants in the distance formula in the following Theorem \hyperref[Theorem 3.11]{3.11}. We also bound $|f_{(\epsilon,(j,\mathbf{k}))}^{(d)}|$ differently, the main motivation for which will become clear in the later Remark \hyperref[Remark 3.17]{3.16}. For now, it can be viewed simply as a more functional analytic bound. Moreover, the following estimate shows it is at least as tight as the previous one:
\begin{align}\frac{1}{a_{14}(\psi,s)}&=\max_{\epsilon\in\alpha}||\psi_{\epsilon}||_{(\dt{C}^s(\mathbb{R}^n))^{\ast}}\notag\\
&=\max_{\epsilon\in\alpha}\sup_{f\in\dt{C}_{1}^{s}(\mathbb{R}^n)}\Biggl|\int_{\mathbb{R}^n}f(\mathbf{x})\psi_{\epsilon}(\mathbf{x})d\mathbf{x}\Biggr|\notag\\
&=\max_{\epsilon\in\alpha}\sup_{f\in\dt{C}_{1}^{s}(\mathbb{R}^n)}\inf_{\mathbf{r}\in\mathbb{R}^n}\Biggl|\int_{\mathbb{R}^n}(f(\mathbf{x})-f(\mathbf{r}))\psi_{\epsilon}(\mathbf{x})d\mathbf{x}\Biggr|\notag\\
&\leqslant\max_{\epsilon\in\alpha}\inf_{\mathbf{r}\in\mathbb{R}^n}\int_{\mathbb{R}^n}|\mathbf{x}-\mathbf{r}|^s|\psi_{\epsilon}(\mathbf{x})|d\mathbf{x}\notag\\
&=\frac{1}{a_{12}(\psi,s)},\notag
\end{align}
where we again use that all $\psi_{\epsilon}$ integrate to 0.\\

Later, we also work with sets similar to $\text{Wav}^d(C_H^s(\mathbb{R}^n))$, derived from the full wavelet
expansion. We conclude this subsection with a similar characterization of $C^s(\mathbb{R}^n)$ in terms of detail coefficients only.\\

\noindent \textbf{Lemma 3.6} \phantomsection\label{Lemma 3.6} (Partial quantitative characterization of $C^s(\mathbb{R}^n)$ via wavelet functions only) Fix $H>0$. For all $f\in C_{H}^s(\mathbb{R}^n)$, when in the full wavelet expansion, the wavelet coefficients are bounded by
\begin{equation}
|f_{(\epsilon,(j,\mathbf{k}))}^{(d)}|\leqslant \frac{H2^{-j(s+\frac{n}{2})}}{a_{14}(\psi,s)} \qquad \epsilon\in\alpha,j\in\mathbb{Z},\mathbf{k}\in\mathbb{Z}^n,
\notag
\end{equation}
where $a_{14}(\psi,s)$ is defined in Lemma \hyperref[Lemma 3.4]{3.4}. Moreover, its wavelet series converges uniformly on compact subsets to $f$.

\noindent \emph{Proof}. The proof of Lemma \hyperref[Lemma 3.4]{3.4} can be carried over without any modification to show the first statement; this is because in showing this direction, besides elementary inequalities, only a change of variables was used, which does not alter when one allows $j<0$ when computing the integrals. For the second statement, see \cite[Chapter 6.4]{bib4}. \qed

\subsection{Original wavelet $s$-Wasserstein distance and its challenges}\label{subsec3.2}
\noindent A key goal of \cite{bib1} is to use Lemma \hyperref[Lemma 3.2]{3.2} to characterize the constraint set in the optimization problem defining the $s$-Wasserstein distance (\ref{eq3.1}) in terms of the size of the wavelet coefficients. In particular, an  immediate consequence of   Lemma \hyperref[Lemma 3.2]{3.2}  is, for all $C_0\geqslant0$ and $C_1>0$, if we define the set
\begin{align}
A_{C_0,C_1}\vcentcolon = \Bigl{\{}(a_{(0,\mathbf{k})})_{\mathbf{k}\in\mathbb{Z}^n}\bigoplus (a_{(\epsilon, (j,\mathbf{k}))})_{\substack{\epsilon \in \alpha, \\j\geqslant0, \mathbf{k}\in\mathbb{Z}^n}}\in & \ \ell^2\ \Big| \ \notag\\
&a_{(0,\mathbf{k})} \leqslant C_0, a_{(\epsilon, (j,\mathbf{k}))} \leqslant C_12^{-j(s+\frac{n}{2})}\Bigr{\}},\notag
\end{align}
then,
\begin{equation}
A_{C_0,C_1}
\subseteq\text{Wav}_{0}^{a}(C_{a_{21}(\phi,s)C_0+a_{22}(\psi,s)C_1}^s(\mathbb{R}^n)).\tag{eq3.4}\label{eq3.4}
\end{equation}
The authors likewise claim that,\footnote{This claim is explicitly stated in the following sentence from their proof of \cite[Theorem 2]{bib1}: ``all functions with $C_H(f)$ less than the lower bound in (14) are included by the constraint.''} for all $C_0\geqslant0$ and $C_1>0$,
\begin{equation}
\text{Wav}_{0}^{a}(C_{a_{12}(\psi,s)C_1}^s(\mathbb{R}^n))\notag\\
\subseteq
A_{C_0,C_1}. \tag{eq3.5}\label{eq3.5}
\end{equation}
While a modified version of this containment is true for certain homogeneous $s$-Hölder spaces and $C_0>0$, as we will describe in detail in Section \ref{subsec3.3}, this claim, as stated, is false, as can be clearly seen in the case $C_0=0$.\footnote{The following example provides a function that belongs to the left-hand side of (\ref{eq3.5}) but does not belong to the right-hand side when $C_0 = 0$: The scaling function $\phi$, being continuously differentiable with compact support, is $s$-Hölder continuous. After multiplying by a $C>0$ sufficiently small, it can have a $s$-Hölder constant less than or equal to $a_{12}(\psi,s)C_1$. However, its approximation coefficients will be strictly positive in view of (\ref{eq2.7}).}

However, if one accepted the containments in (\ref{eq3.4}, \ref{eq3.5}), it is clear that the $s$-Wasserstein distance between $\mu$ and $\nu$ could be bounded from above and below in terms of constant multiples of the following optimization problem:
$$\sup_{a\in A_{C_0,C_1}} \Biggl( \sum_{\mathbf{k} \in \mathbb{Z}^n} a_{(0,\mathbf{k})} (p_{\mu}-p_{\nu})_{(0,\mathbf{k})}^{(a)} + \sum_{\substack{\epsilon \in \alpha,\\ j \geqslant 0, \mathbf{k} \in \mathbb{Z}^n}} a_{(\epsilon, (j, \mathbf{k}))} (p_{\mu}-p_{\nu})_{(\epsilon, (j, \mathbf{k}))}^{(d)} \Biggr).$$
This optimization problem can be easily solved to obtain the solution 
\begin{equation}
\sum_{\mathbf{k} \in \mathbb{Z}^n} C_0|(p_{\mu}-p_{\nu})_{(0,\mathbf{k})}^{(a)}| + \sum_{\substack{\epsilon\in\alpha,\\ j\geqslant 0,\  \mathbf{k}\in\mathbb{Z}^n}}C_12^{-j(s+\frac{n}{2})}|(p_{\mu}-p_{\nu})_{(\epsilon, (j, \mathbf{k}))}^{(d)}|,\notag
\end{equation}
where the square summability of the sequence of coefficients that attains the supremum is ensured by the assumption that $\mu$, $\nu$, and the wavelet system $(\phi, \psi)$ are all compactly supported. This closed form expression provides the rationale behind the main theorem in \cite{bib1}. However, it relies on the incorrect containment in (\ref{eq3.5}). We now recall the statement of theorem from \cite{bib1} and provide an example illustrating that, not only is the rationale incorrect, as explained above, but the theorem statement itself is also false.\\

\noindent \textbf{Theorem 3.7}\phantomsection\label{Theorem 3.7} (Original wavelet $s$-Wasserstein distance for $0<s<1$\footnote{Our presentation of the theorem differs from \cite[Theorem 2]{bib1} in a few aspects. First, in the original version, the measures are not required to have square-integrable densities, nor does the system $(\phi,\psi)$ need to be compactly supported; we introduce these hypotheses to ensure other parts of the logic flow properly. Second, there is a result for discrete measures in the case of $s=1$; see Remark \hyperref[Remark 3.9]{3.9}.} \cite[Theorem 2]{bib1}) Fix $0<s<1$ and a compactly supported $r$-regular orthonormal wavelet system $(\phi,\psi)$ with $r\geqslant1$. Given $d\mu = p_{\mu}(\mathbf{x})d\mathbf{x}$, $d\nu = p_{\nu}(\mathbf{x})d\mathbf{x}\in\mathcal{P}_{c,\text{ab2}}(\mathbb{R}^n)$, define the \textbf{original wavelet $\pmb{s}$-Wasserstein distance} $(W_{\text{wav,ori}}^{(\phi,\psi,C_0,C_1)})_{s}(\mu,\nu)$ \textbf{with} $\pmb{C_0\geqslant0}$ \textbf{and} $\pmb{C_1>0}$ as
\begin{equation}
(W_{\text{wav,ori}}^{(\phi,\psi,C_0,C_1)})_{s}(\mu,\nu) \vcentcolon  = \sum_{\mathbf{k} \in \mathbb{Z}^n} C_0|(p_{\mu}-p_{\nu})_{(0,\mathbf{k})}^{(a)}| + \sum_{\substack{\epsilon\in\alpha,\\ j\geqslant 0,\  \mathbf{k}\in\mathbb{Z}^n}}C_12^{-j(s+\frac{n}{2})}|(p_{\mu}-p_{\nu})_{(\epsilon, (j, \mathbf{k}))}^{(d)}|.\notag
\end{equation}
Then $(W_{\text{wav,ori}}^{(\phi,\psi,C_0,C_1)})_{s}$ is equivalent to the classical $s$-Wasserstein distance $W_s$ on $\mathcal{P}_{c,\text{ab2}}(\mathbb{R}^n)$, in the sense that
\begin{equation}
C_{\psi,s,C_1}W_s(\mu, \nu)\leqslant (W_{\text{wav,ori}}^{(\phi,\psi,C_0,C_1)})_{s}(\mu, \nu) \leqslant \tilde{C}_{\phi,\psi,s,C_0,C_1}W_s(\mu, \nu)\tag{eq3.6}\label{eq3.6}
\end{equation}
for all $\mu$, $\nu\in\mathcal{P}_{c,\text{ab2}}(\mathbb{R}^n)$, where $C_{\psi,s,C_1} \vcentcolon= a_{12}(\psi,s)C_1$, and $\tilde{C}_{\phi,\psi,s,C_0,C_1} \vcentcolon = a_{21}(\phi,s)C_0+$ $a_{22}(\psi,s)C_1$. When $C_0=0$ and $C_1=1$, $(W_{\text{wav,ori}}^{(\phi,\psi,C_0,C_1)})_{s}$ is simply called the \textbf{original wavelet $\pmb{s}$-Wasserstein distance} and is abbreviatd as $(W_{\text{wav,ori}}^{\psi})_{s}$.\\

\noindent The authors in \cite{bib1} particularly advocate the choice $C_0 = 0$ ``because this gives the tightest bounds'' and subsequently also set $C_1 = 1$. We now provide a counterexample to the theorem in the case. Below, we shall use the following notation for the translation operator on $\mathcal{P}_{c,\text{ab2}}(\mathbb{R}^n)$:
$$\tau_{\mathbf{a}}(\mu) \vcentcolon  = p_{\mu}(\mathbf{x}-\mathbf{a})d\mathbf{x}$$
for $\mathbf{a}\in\mathbb{R}^n$, $d\mu = p_{\mu}(\mathbf{x})d\mathbf{x}\in\mathcal{P}_{c,\text{ab2}}(\mathbb{R}^n)$. \\

\noindent \textsc{Example 3.8}\phantomsection\label{Example 3.8} (Counterexample to Theorem \hyperref[Theorem 3.7]{3.7}) Assume $n=1$. Let our system $(\phi,\psi)$ satisfy $\text{supp}(\psi)\subseteq [0,4]$, which is satisfied, for example, by the Daubechies 4 wavelet system, described in Daubechies \cite[Chapters 6 \& 7]{bib7}. Let $d\mu = p_\mu(x)dx \vcentcolon  = \mathbbm{1}_{[0,1]} dx\in\mathcal{P}_{c,\text{ab2}}(\mathbb{R})$ be the Lebesgue measure restricted to [0,1]. Then, for all integers $m\geqslant0$, $j\geqslant0$, and $k$, we observe that
$$\left|\int_{\mathbb{R}}(p_{\mu}-p_{\tau_{4+m}(\mu)})(x)\psi_{(j,k)}^{\ast}(x)dx\right|$$
equals either $$\left|\int_{\mathbb{R}}p_{\mu}(x)\psi_{(j,k)}^{\ast}(x)dx\right| \qquad \text{or} \qquad \left|\int_{\mathbb{R}}p_{\tau_{4+m}(\mu)}(x)\psi_{(j,k)}^{\ast}(x)dx\right|,$$
and the other must be zero. Consequently, for all integers $m\geqslant0$,
\begin{align}
    (W_{\text{wav,ori}}^{\psi})_s(\mu,\tau_4 (\mu)) &= \sum_{j\geqslant0,k\in\mathbb{Z}}2^{-j(s+\frac{1}{2})}|(p_{\mu}-p_{\tau_4(\mu)})_{(j, k)}^{(d)}|\notag\\
    &= \sum_{j\geqslant0,k\in\mathbb{Z}}2^{-j(s+\frac{1}{2})}(|(p_{\mu})_{(j, k)}^{(d)}|+|(p_{\tau_4(\mu)})_{(j, k)}^{(d)}|)\notag\\
    &= \sum_{j\geqslant0,k\in\mathbb{Z}}2^{-j(s+\frac{1}{2})}(|(p_{\mu})_{(j, k)}^{(d)}|+|(p_{\tau_{4+m}(\mu)})_{(j, k)}^{(d)}|)\notag\\
    &=\sum_{j\geqslant0,k\in\mathbb{Z}}2^{-j(s+\frac{1}{2})}|(p_{\mu}-p_{\tau_{4+m}(\mu)})_{(j, k)}^{(d)}|\notag\\
    &=(W_{\text{wav,ori}}^{\psi})_s(\mu,\tau_{4+m} (\mu))\tag{eq3.7}\label{eq3.7}
\end{align}
where the third equality can be shown by using a change of variables on (\ref{eq2.8}), since $(p_{\tau_4(\mu)})_{(j,k)}^{(d)} = (p_{\tau_{4+m}(\mu)})_{(j,2^jm + k)}^{(d)}$ and $2^jm + k\in\mathbb{Z}$. 

By equation (\ref{eq2.1}), $W_s(\mu,\tau_y (\mu))\to\infty$ as $y\to\infty$; therefore, the lower bound in (\ref{eq3.6}) implies that the right hand side of (\ref{eq3.7}) should diverge to $\infty$ as $m\to\infty$. However, the upper bound in (\ref{eq3.6}) implies that the left-hand side of (\ref{eq3.7}) is bounded independent of $m$. This is a contradiction. For other $n$ and wavelet systems, the notations are slightly more complicated but the same argument applies.\\

\noindent Two more remarks follow.\\

\noindent \textsc{Remark 3.9}\phantomsection\label{Remark 3.9} ($(W_{\text{wav,ori}}^{(\phi,\psi,C_0,C_1)})_{s}$ for $s=1$) Due to Remark \hyperref[Remark 3.3]{3.3}, it is clear that, even setting aside the logical gaps of Theorem \hyperref[Theorem 3.7]{3.7},  the current technique in establishing a wavelet equivalent $s$-Wasserstein distance can not extend to the case of $s=1$. The authors in \cite{bib1} acknowledged a similar issue, yet asserted that the statement in Theorem \hyperref[Theorem 3.7]{3.7} remains valid for arbitrary discrete probability measures when $s=1$; see \cite[Theorem 2]{bib1}. After examination, we think that they do not completely justify their final conclusion, unless they meant to restrict to discrete measures with support points on a grid. Given that this discussion veers significantly into abstract harmonic analysis and diverges from our main focus on an absolutely continuous measure framework, we chose not to pursue this further and do not include the $s=1$ case in Theorem \hyperref[Theorem 3.7]{3.7}, at least theoretically.\footnote{Our numerical simulations do include $s=1$, see Section \ref{sec4}.}\\

\noindent \textsc{Remark 3.10}\phantomsection\label{Remark 3.10} (Simulations of $(W_{\text{wav,ori}}^{\psi})_{s}$ and intuitive interpretation) As a complement to the preceding discussion of the theoretical difficulties of $(W_{\text{wav,ori}}^{\psi})_{s}$, we also examine its properties numerically in Section \ref{subsec4.2}. Since $(W_{\text{wav,ori}}^{\psi})_{s}$ does not precisely converge to $W_s$ in any sense, we do not expect them to exactly agree. However, a reasonable baseline to assess whether $(W_{\text{wav,ori}}^{\psi})_{s}$ behaves numerically as suitable proxy for $W_s$ is to expect that the ratio $(W_{\text{wav,ori}}^{\psi})_s/W_s$ is roughly a constant between all measures in every ``simple'' data set, each characterized by a specific pattern. If so, after normalization, the wavelet distance would capture the classical distance well, at least in the simple cases. Thus, our simulations were done for translations and dilations of common measures, the two most elementary types of measure transport, and the results were compared with the classical $W_s$, after applying normalizing constants. 

Unfortunately, the theoretical difficulties of $(W_{\text{wav,ori}}^{\psi})_s$ described earlier in this section are accompanied by poor numerical performance. In each fixed data set, the ratio $(W_{\text{wav,ori}}^{\psi})_s/W_s$ decreases rapidly as the support of $\mu-\nu$ becomes large. This is not surprising since, for any $d\mu = p_{\mu}(\mathbf{x})d\mathbf{x}$, $d\nu = p_{\nu}(\mathbf{x})d\mathbf{x}$, the original wavelet $s$-Wasserstein distance seeks to use
\begin{equation}\sum_{\substack{\epsilon\in\alpha,\\ j\geqslant 0,\  \mathbf{k}\in\mathbb{Z}^n}}2^{-j(s+\frac{n}{2})}\text{sgn}((p_{\mu}-p_{\nu})_{(\epsilon, (j, \mathbf{k}))}^{(d)})\psi_{(\epsilon, (j, \mathbf{k}))},\tag{eq3.8}\label{eq3.8}\end{equation}
to approximate the $f_{\text{max}}$ that attains the supremum in (\ref{eq3.1}) for some $H>0$. Nevertheless, without the contribution from either the level 0 approximation coefficients or the negative level detail coefficients, the final constructed (\ref{eq3.8}) may miss components to approximate $f_{\text{max}}$ well for any $H>0$, and this effect becomes severe when the domain of measures' support is large. The counterexample constructed in Example \hyperref[Example 3.8]{3.8} is exploiting this point rigorously.

\subsection{An attempted reformulation}\label{subsec3.3}
We now show that if one instead requires $C_0>0$ in the above reasoning, then the previous result can be amended to some extent: $(W_{\text{wav,ori}}^{(\phi,\psi,C_0,C_1)})_{s}$ will be an equivalent distance to  $W_s$, for $\mathcal{P}_{\text{ab2}}(K)$ on compact domains $K\subseteq\mathbb{R}^n$.

The argument to show it relies on three observations: first, in (\ref{eq3.1}), the values of $f$ outside of any $K\supseteq \text{supp}(p_{\mu} - p_{\nu})$ can be freely altered without changing the value of integral; second, every function in $C_H^s(K)$ is a restriction of a function in $C_H^s(\mathbb{R}^n)$ with the same supremum and $s$-Hölder constant, due to the Kirszbraun theorem; third, since the integral of $p_{\mu}-p_{\nu}$ vanishes, modifying $f$ on $K$ by addition of a constant does not change the value of the integral in (\ref{eq3.1}).

Motivated by these considerations, we will consider the homogeneous space $\dt{C}^{s}_{H}(K)$, instead of the original H\"older space $C^s_H(K)$. Furthermore, we commit a mild abuse of notation, allowing $\dt{C}^{s}_{H}(K)$ to denote the realizations via functions that satisfy $f(\mathbf{0}) = 0$ and are understood to be continuously extended from $K$ to $\mathbb{R}^n$, while preserving both the supremum on $K$ and the $s$-Hölder constant. With these conventions, the dual formulation of the $s$-Wasserstein distance (\ref{eq3.1}) can be rewritten as
\begin{equation}
HW_s(\mu, \nu) = \sup_{f\in \dt{C}^s_{H}(K)\cap L^2(\mathbb{R}^n)}\int_{K}f(\mathbf{x}) (p_{\mu}(\mathbf{x})-p_{\nu}(\mathbf{x}))d\mathbf{x}.\tag{eq3.9}\label{eq3.9}\end{equation}

Together with Lemma \hyperref[Lemma 3.4]{3.4}, this alternative formulation allows us to repair the problems of (\ref{eq3.5}) and Theorem \hyperref[Theorem 3.7]{3.7}, as we now describe.\\

\noindent \textbf{Theorem 3.11}\phantomsection\label{Theorem 3.11} (Alternative formulation of the original wavelet $s$-Wasserstein distance for $0<s<1$) Fix a compact set $K\subseteq\mathbb{R}^n$, $0<s<1$, and a compactly supported $r$-regular orthonormal wavelet system $(\phi,\psi)$ with $r\geqslant1$; let $B\vcentcolon = \text{diam}(K)$. Then, given $d\mu = p_{\mu}(\mathbf{x})d\mathbf{x}$, $d\nu = p_{\nu}(\mathbf{x})d\mathbf{x}\in\mathcal{P}_{\text{ab2}}(K)$, for any $C_0>0$ and $C_1>0$, $(W_{\text{wav,ori}}^{(\phi,\psi,C_0,C_1)})_{s}$ is equivalent to the classical $s$-Wasserstein distance $W_s$ on $\mathcal{P}_{\text{ab2}}(K)$. In particular, there exist positive constants $D_{\phi,\psi,s,C_0,C_1,B}$, $\tilde{D}_{\phi,\psi,s,C_0,C_1}$
\begin{align}
D_{\phi,\psi,s,C_0,C_1,B} &\vcentcolon= \min{(a_{13}(\phi)C_0/B^s  ,a_{14}(\psi,s)C_1)}\notag\\
\tilde{D}_{\phi,\psi,s,C_0,C_1} &\vcentcolon= a_{21}(\phi,s)C_0+a_{22}(\psi,s)C_1\notag
\end{align}
so that
\begin{equation}
D_{\phi,\psi,s,C_0,C_1,B}W_s(\mu, \nu)\leqslant (W_{\text{wav,ori}}^{(\phi,\psi,C_0,C_1)})_s(\mu, \nu) \leqslant \tilde{D}_{\phi,\psi,s,C_0,C_1}W_s(\mu, \nu)\tag{eq3.10}\label{eq3.10}
\end{equation}
for all $\mu$, $\nu\in\mathcal{P}_{\text{ab2}}(K)$.

\noindent \emph{Proof}. By Remark \hyperref[Remark 2.3]{2.3}, we may assume, without loss of generality, that $K$ contains $\mathbf{0}$. Note that, for every $f\in C_{H}^s(K)$ such that $f(\mathbf{0}) = 0$, its supremum on $K$ is bounded by $HB^s$. If we denote
\begin{equation}
\text{Wav}_{j_0}^{a}(\dt{C}^s_{H}(K)) \vcentcolon  = \Bigl{\{}(f_{(j_0,\mathbf{k})}^{(a)})_{\mathbf{k}\in\mathbb{Z}^n}\bigoplus (f_{(\epsilon, (j,\mathbf{k}))}^{(d)})_{\substack{\epsilon \in \alpha,\\ j \geqslant j_0, \mathbf{k} \in \mathbb{Z}^n}}\in \ell^2\ \Big|\ f \in \dt{C}^s_{H}(K) \Bigr{\}},\notag
\end{equation}
with a realized $\dt{C}^s_{H}(K)$, then by Lemma \hyperref[Lemma 3.2]{3.2} and Lemma \hyperref[Lemma 3.4]{3.4},
\begin{align}
\text{Wav}_{0}^{a}(\dt{C}_{D_{\phi,\psi,s,C_0,C_1,B}}^s(K))\subseteq &\Bigl{\{}(a_{(0,\mathbf{k})})_{\mathbf{k}\in\mathbb{Z}^n}\bigoplus (a_{(\epsilon, (j,\mathbf{k}))})_{\substack{\epsilon \in \alpha,\\ j\geqslant0, \mathbf{k}\in\mathbb{Z}^n}}\in \ell^2\ \Big| \notag\\& \ \ a_{(0,\mathbf{k})} \leqslant C_0, \ a_{(\epsilon, (j,\mathbf{k}))} \leqslant C_12^{-j(s+\frac{n}{2})}\Bigr{\}}\notag\\
\subseteq&\text{Wav}_{0}^{a}(C_{\tilde{D}_{\phi,\psi,s,C_0,C_1}}^s(\mathbb{R}^n)).\notag
\end{align}
Using (\ref{eq3.1}), (\ref{eq3.9}), and following a similar logic in the final steps in obtaining the original wavelet $s$-Wasserstein distance as described before Theorem \hyperref[Theorem 3.7]{3.7}, we arrive at (\ref{eq3.10}).\qed\\

\noindent We now comment on its numerical performance.\\

\noindent \textsc{Remark 3.12}\phantomsection\label{Remark 3.12} (Simulations of $(W_{\text{wav,ori}}^{(\phi,\psi,C_0,C_1)})_s$ and intuitive interpretation) In Section \ref{subsec4.2}, we conduct analogous simulations to those described in Remark \hyperref[Remark 3.10]{3.10}, but now with $C_0>0$ sufficiently large, with respect to the support of the measures and $C_1=1$. The results are slightly better than the original ones but still unsatisfactory; the ratio $(W_{\text{wav,ori}}^{(\phi,\psi,C_0,C_1)})_s/W_s$ is not roughly constant between all measures within each fixed data set. 

We believe this is for two reasons. First, this time, the formulation of the wavelet $s$-Wasserstein distance approximates the $f_{\text{max}}$ that achieves the supremum in (\ref{eq3.1}) for some $H>0$ by
\begin{equation}\sum_{\mathbf{k}\in\mathbb{Z}^n}C_0\text{sgn}((p_{\mu}-p_{\nu})_{(0,\mathbf{k})}^{(a)})\phi_{(0,\mathbf{k})}+\sum_{\substack{\epsilon\in\alpha,\\ j\geqslant 0,\  \mathbf{k}\in\mathbb{Z}^n}}C_12^{-j(s+\frac{n}{2})}\text{sgn}((p_{\mu}-p_{\nu})_{(\epsilon, (j, \mathbf{k}))}^{(d)})\psi_{(\epsilon, (j, \mathbf{k}))},\tag{eq3.11}\label{eq3.11}
\end{equation}
which strongly depends on the artificial choice of the combination of constants $C_0$ and $C_1$ that provides the scale for the approximation and detail coefficients of this constructed (\ref{eq3.11}). This is unnatural from the perspective of the dual formulation of the $s$-Wasserstein distance in (\ref{eq3.1}): 
the approximation coefficients characterize the inhomogeneous part of a function (size) and the detail coefficients characterize the homogeneous part of a function (oscillation). A choice on the combination of $C_0$ and $C_1$ imposes an extra constraint on the ratio of size and oscillation for the function (\ref{eq3.11}). However, the constraint in defining the $W_s$ is imposed only on the homogeneous part of $f$ and should also scale linearly with such part.

A second potential reason for the poor numerical performance is, using the notation in Example \hyperref[Example 3.8]{3.8} and mimicking the proof, one can easily show that  
\begin{equation}(W_{\text{wav,ori}}^{(\phi,\psi,C_0,C_1)})_s(\mu,\tau_4 (\mu)) =  (W_{\text{wav,ori}}^{(\phi,\psi,C_0,C_1)})_s(\mu,\tau_{4+m} (\mu))\tag{eq3.12}\label{eq3.12}
\end{equation}for all integers $m\geqslant0$, after one chooses a large enough domain $K$ depending on $m$. This periodic phenomenon remains even after we reestablish the equivalence on the distance level; nevertheless, it is not how $W_s$ behaves: it should be strictly increasing with $m$. The proof relies strongly on the fact that there is a maximum diameter for the supports of all the scaling functions and wavelets with level $j\geqslant0$. This shows that, without less localized scaling or wavelet functions, the wavelet $s$-Wasserstein distance is inherently incompatible with translation, an essential feature of $W_s$, when the domain is large enough.

\subsection{New wavelet $s$-Wasserstein distance}\label{subsec3.4}

Given the challenges of the original wavelet $s$-Wasserstein distance discussed in Section \ref{subsec3.3}, we now introduce a final new formulation of the wavelet $s$-Wasserstein distance, starting from the full wavelet expansion of the Hölder continuous functions via detail coefficients only. This removes the need to choose the artificial combination of constant $C_0$ and $C_1$ in our alternative formulation of the original distance in Section \ref{subsec3.3}; moreover, by including detail coefficients from less localized wavelets, it prevents periodic behavior under translations that occurred in the previous formulations.\\

\noindent \textbf{Theorem 3.13}\phantomsection\label{Theorem 3.13} (New wavelet $s$-Wasserstein distance for $0<s<1$) Fix a compact set $K\subseteq \mathbb{R}^n$, $0<s<1$, and a compactly supported $r$-regular orthonormal wavelet system $(\phi,\psi)$ with $r\geqslant1$. Given $d\mu = p_{\mu}(\mathbf{x})d\mathbf{x}$, $d\nu = p_{\nu}(\mathbf{x})d\mathbf{x}\in\mathcal{P}_{\text{ab2}}(K)$, define the \textbf{new wavelet Wasserstein $\pmb{s}$-distance} $(W_{\text{wav}}^{(\psi,j_0)})_s(\mu,\nu)$ \textbf{with the lowest level $\pmb{j_0\in\mathbb{Z}}$} between them as
\begin{equation}
(W_{\text{wav}}^{(\psi,j_0)})_s(\mu,\nu) \vcentcolon  = \sum_{\substack{\epsilon\in\alpha,\\ j\geqslant j_0,\  \mathbf{k}\in\mathbb{Z}^n}}2^{-j(s+\frac{n}{2})}|(p_{\mu}-p_{\nu})_{(\epsilon, (j, \mathbf{k}))}^{(d)}|.\notag
\end{equation}
Then, for all $\eta>0$, there exists sufficiently small integer $j_0\in\mathbb{Z}$, depending on the domain $K$, $s$, and the system $(\phi,\psi)$, such that
\begin{equation}
E_{\psi,s}W_s(\mu, \nu)-\eta\leqslant (W_{\text{wav}}^{(\psi,j_0)})_s(\mu, \nu) \leqslant \tilde{E}_{\psi,s,j_0}W_s(\mu, \nu)\tag{eq3.13}\label{eq3.13}
\end{equation}
where 
\begin{align}E_{\psi,s} &\vcentcolon= a_{14}(\psi,s),\notag\\
\tilde{E}_{\psi,s,j_0} &\vcentcolon= 
\begin{cases}
a_{22}(\psi,s) + |j_0|a_{23}(\psi,s) & \text{if} \ j_0<0 ,\notag\\
a_{22}(\psi,s) & \text{if} \ j_0\geqslant0.
\end{cases}
\end{align}
The constants $a_{14}(\psi,s)$ and $a_{22}(\psi,s)$ are defined as before, while
\begin{equation}
a_{23}(\psi,s)\vcentcolon  = ||\textstyle\sum\psi^s||_{\infty}\vcentcolon  = \sup_{\substack{\mathbf{x}\neq\mathbf{y}\\\eta_k\in\{+1,-1\}}}\left|\frac{\sum_{\mathbf{k}\in\mathbb{Z}^n,\epsilon\in\alpha}\eta_k\psi_{(\epsilon,(0,\mathbf{k}))}(\mathbf{x})-\sum_{\mathbf{k}\in\mathbb{Z}^n,\epsilon\in\alpha}\eta_k\psi_{(\epsilon,(0,\mathbf{k}))}(\mathbf{y})}{|\mathbf{x}-\mathbf{y}|^s} \right|.\notag
\end{equation}
\noindent Besides the alternative formulation (\ref{eq3.9}) in Section \ref{subsec3.2} and the preparatory Lemma \hyperref[Lemma 3.6]{3.6} displayed in Section \ref{subsec3.1}, our proof of this theorem relies on the following important lemma.\\

\noindent \textbf{Lemma 3.14}\phantomsection\label{Lemma 3.14} Given $\eta>0$ and a compact set $K$, there exists sufficiently small $j_0\in\mathbb{Z}$, depending on the domain $K$, $s$, and wavelet system $(\phi,\psi)$, such that for every $f\in[f]\in \dt{C}_{H}^s(\mathbb{R}^n)$, with the space interpreted as the traditional quotient space, defining
$$g(\mathbf{x}) \vcentcolon  =  \sum_{\substack{\epsilon\in\alpha,\\ j\geqslant j_0,\  \mathbf{k}\in\mathbb{Z}^n}}f_{(\epsilon, (j, \mathbf{k}))}^{(d)}\psi_{(\epsilon, (j, \mathbf{k}))}(\mathbf{x}),$$
we have $||g-\tilde{f}||_{L^{\infty}(K)}<\eta$ for some $\tilde{f}\in[f]$.

\noindent \emph{Proof}. Denote 
$$f_j(\mathbf{x}) \vcentcolon  =  \sum_{\epsilon\in\alpha, \  \mathbf{k}\in\mathbb{Z}^n}f_{(\epsilon, (j, \mathbf{k}))}^{(d)}\psi_{(\epsilon, (j, \mathbf{k}))}(\mathbf{x})$$
and $B\vcentcolon  = \text{diam}(K)$. Fix $\mathbf{x_0}\in K$. Then for $\mathbf{x}\in K$,
\begin{align}
|f_{j}(\mathbf{x})-f_{j}(\mathbf{x_0})| &= \left|\sum_{\epsilon\in\alpha, \  \mathbf{k}\in\mathbb{Z}^n}f_{(\epsilon, (j, \mathbf{k}))}^{(d)}\psi_{(\epsilon, (j, \mathbf{k}))}(\mathbf{x}) - \sum_{\epsilon\in\alpha, \  \mathbf{k}\in\mathbb{Z}^n}f_{(\epsilon, (j, \mathbf{k}))}^{(d)}\psi_{(\epsilon, (j, \mathbf{k}))}(\mathbf{x_0})\right| \notag \\
& = \left|2^{\frac{nj}{2}}\sum_{\epsilon\in\alpha, \  \mathbf{k}\in\mathbb{Z}^n}f_{(\epsilon, (j, \mathbf{k}))}^{(d)}(\psi_{\epsilon}(2^j\mathbf{x}-\mathbf{k}) - \psi_{\epsilon}(2^j\mathbf{x_0}-\mathbf{k}))\right|. \notag
\end{align}
Since $\{\psi_{\epsilon}\}_{\epsilon\in\alpha}$ are at least 1-regular and compactly supported, they are 1-Hölder with a maximum Hölder constant bound, say $H_{\psi}$. Then, together with Lemma \hyperref[Lemma 3.6]{3.6},
\begin{align}
|f_{j}(\mathbf{x})-f_{j}(\mathbf{x_0})| \leqslant &2^{\frac{nj}{2}}\frac{H2^{-j(s+\frac{n}{2})}}{a_{14}(\psi,s)}\sum_{\epsilon\in\alpha, \  \mathbf{k}\in\mathbb{Z}^n}|\psi_{\epsilon}(2^j\mathbf{x}-\mathbf{k}) - \psi_{\epsilon}(2^j\mathbf{x_0}-\mathbf{k}))| \notag\\
\leqslant & \frac{H2^{-j(s-1)}}{a_{14}(\psi,s)}CBH_{\psi} \tag{eq3.14}\label{eq3.14}
\end{align}
for some constant $C$ independent of $j$ since $\{\psi_{\epsilon}\}_{\epsilon\in\alpha}$ are compactly supported (which results in the right-hand side of the first inequality in (\ref{eq3.14}) to have at most finitely many nonzero terms for any $\mathbf{x} \in K$, where the maximum number of terms is bounded uniformly in $\mathbf{x}$). As $2^{s-1}<1$, the series 
$$\sum_{j\leqslant j'}\frac{H2^{-j(s-1)}}{a_{14}(\psi,s)}CBH_{\psi}$$
converges for all $j'\in\mathbb{Z}$. Then there exists a $j_0\in\mathbb{Z}$ such that 
$$\Biggl|\sum_{j< j_0}f_{j}(\mathbf{x})-f_{j}(\mathbf{x_0})\Biggr|\leqslant\sum_{j< j_0}\frac{H2^{-j(s-1)}}{a_{14}(\psi,s)}CBH_{\psi}<\eta.$$
Now, let 
$\tilde{f}(\mathbf{x}) = f(\mathbf{x})-\sum_{j< j_0} f_j(\mathbf{x_0})$, then we have $||g-\tilde{f}||_{L^{\infty}(K)}<\eta$ by the above estimation and $\tilde{f}\in[f]$ since we only subtract $f$ by a constant to get $\tilde{f}$.\qed\\

\noindent We now prove Theorem \hyperref[Theorem 3.13]{3.13}.

\noindent \emph{Proof of Theorem \hyperref[Theorem 3.13]{3.13}}.
Note that, given two functions $\tilde{f}$ and $g$ with $||\tilde{f}-g||_{L^{\infty}(K)}<\eta'$ for some $\eta'>0$, we have
\begin{align}
&\left|\int_{\mathbb{R}^n}\tilde{f}(\mathbf{x}) (p_{\mu}(\mathbf{x})-p_{\nu}(\mathbf{x}))d\mathbf{x} - \int_{\mathbb{R}^n}g(\mathbf{x}) (p_{\mu}(\mathbf{x})-p_{\nu}(\mathbf{x}))d\mathbf{x}\right|\notag \\
\leqslant &\int_{\mathbb{R}^n}|\tilde{f}(\mathbf{x})-g(\mathbf{x})| |p_{\mu}(\mathbf{x})-p_{\nu}(\mathbf{x})|d\mathbf{x} \notag\\
\leqslant &2\eta'.\label{eq3.15}\tag{eq3.15}
\end{align}
for all $p_{\mu}(\mathbf{x})d\mathbf{x}, p_{\nu}(\mathbf{x})d\mathbf{x}\in\mathcal{P}_{\text{ab2}}(K)$. Fix $\eta>0$. By the above estimate and Lemma \hyperref[Lemma 3.14]{3.14}, we may replace every $f$ in (\ref{eq3.9}) by its full wavelet expansion, but only summing above a sufficiently small index $j_0$, denoted as $g$, and the value of the new integral will differ from the original one by at most $\eta$. Since this estimate is uniform in $f$, the difference of the final value with the actual $W_s(\mu,\nu)$ is also bounded by $\eta$.

Next, observe that, if we denote 
\begin{equation}
\text{Wav}_{j_0}^{d}(\dt{C}_H^s(K)) \vcentcolon  = \left\{ (f_{(\epsilon, (j,\mathbf{k}))}^{(d)})_{\substack{\epsilon\in\alpha,\\j\geqslant j_0,\mathbf{k}\in\mathbb{Z}^n}}\in \ell^2\ \Big|\ f \in \dt{C}^{s}_{H}(K) \right\},\notag
\end{equation}
with a realized $\dt{C}_H^s(K)$, then, by Lemma \hyperref[Lemma 3.6]{3.6}, we have 
\begin{equation}
\text{Wav}_{j_0}^{d}(\dt{C}_{E_{\psi,s}}^s(K))\subseteq \Bigl{\{}(a_{(\epsilon, (j,\mathbf{k}))})_{\epsilon \in \alpha, j\geqslant j_0, \mathbf{k}\in\mathbb{Z}^n}\in \ell^2\ \Big|  \ a_{(\epsilon, (j,\mathbf{k}))} \leqslant 2^{-j(s+\frac{n}{2})}\Bigr{\}}.\notag
\end{equation}
Combining the two points above and following a similar logic as described before Theorem \hyperref[Theorem 3.7]{3.7}, we have the first inequality in (\ref{eq3.13}).

For the other inequality, given any $p_{\mu}(\mathbf{x})d\mathbf{x}, p_{\nu}(\mathbf{x})d\mathbf{x}\in\mathcal{P}_{\text{ab2}}(K)$, consider the following function, which we interpret as the function used to approximate the $f_{\text{max}}$ that achieves the supremum in (\ref{eq3.1}) for some $H$,
\begin{align}
&\sum_{\substack{\epsilon\in\alpha,\\ j\geqslant j_0,\  \mathbf{k}\in\mathbb{Z}^n}}2^{-j(s+\frac{n}{2})}\text{sgn}((p_{\mu}-p_{\nu})_{(\epsilon, (j, \mathbf{k}))}^{(d)})\psi(\mathbf{x})_{(\epsilon, (j, \mathbf{k}))}\notag \\ 
=&\sum_{\substack{\epsilon\in\alpha,\\ j\geqslant 0,\  \mathbf{k}\in\mathbb{Z}^n}}2^{-j(s+\frac{n}{2})}\text{sgn}((p_{\mu}-p_{\nu})_{(\epsilon, (j, \mathbf{k}))}^{(d)})\psi(\mathbf{x})_{(\epsilon, (j, \mathbf{k}))} + \notag\\
&\sum_{\substack{\epsilon\in\alpha,\\ j_0\leqslant j <0,\\  \mathbf{k}\in\mathbb{Z}^n}}2^{-j(s+\frac{n}{2})}\text{sgn}((p_{\mu}-p_{\nu})_{(\epsilon, (j, \mathbf{k}))}^{(d)}\psi(\mathbf{x})_{(\epsilon, (j, \mathbf{k}))}.\tag{eq3.16}\label{eq3.16}
\end{align}
By Lemma \hyperref[Lemma 3.2]{3.2} and the fact that both the wavelet system and $p_{\mu}-p_{\nu}$ are compactly supported, the first term on the right-hand side of (\ref{eq3.16}) belongs to $C_{a_{22}(\psi,s)}^s(\mathbb{R}^n)\cap L^2(\mathbb{R}^n)$. For each summand in the second term, by definition of our $a_{23}(\psi,s)$,
\begin{align}
&\Biggl|\sum_{\mathbf{k}\in\mathbb{R}^n,\epsilon\in\alpha}2^{-j(s+\frac{n}{2})}\text{sgn}((p_{\mu}-p_{\nu})_{(\epsilon, (j, \mathbf{k}))}^{(d)})\psi(\mathbf{x})_{(\epsilon, (j, \mathbf{k}))} - \notag\\&\ \sum_{\mathbf{k}\in\mathbb{R}^n,\epsilon\in\alpha}2^{-j(s+\frac{n}{2})}\text{sgn}((p_{\mu}-p_{\nu})_{(\epsilon, (j, \mathbf{k}))}^{(d)})\psi(\mathbf{y})_{(\epsilon, (j, \mathbf{k}))}\Biggr|\notag\\
=&\Biggl|\sum_{\mathbf{k}\in\mathbb{R}^n,\epsilon\in\alpha}2^{-j(s+\frac{n}{2})}\text{sgn}((p_{\mu}-p_{\nu})_{(\epsilon, (j, \mathbf{k}))}^{(d)})2^{\frac{nj}{2}}(\psi_{\epsilon}(2^j\mathbf{x}-\mathbf{k})-\psi_{\epsilon}(2^j\mathbf{y}-\mathbf{k}))\Biggr|\notag\\
\leqslant&2^{-js}a_{23}(\psi,s)|(2^j\mathbf{x}-\mathbf{k})-(2^j\mathbf{y}-\mathbf{k})|^s\notag\\
=&a_{23}(\psi,s)|\mathbf{x}-\mathbf{y}|^s.\notag
\end{align}
Together again with the fact that both the wavelet system and $p_{\mu}-p_{\nu}$ are compactly supported, we see that each term in the second sum on the right-hand side of (\ref{eq3.16}) belongs to $C_{a_{23}(\psi,s)}^s(\mathbb{R}^n)\cap L^2(\mathbb{R}^n)$. Since the Hölder constant is a seminorm on the Hölder space $C^s(\mathbb{R}^n)$, the overall expression in (\ref{eq3.16}) belongs to $C_{\tilde{E}_{\psi,s,j_0}}^s(\mathbb{R}^n)\cap L^2(\mathbb{R}^n)$; then the second inequality follows.\qed\\

We now make some important remarks and describe a few open questions, which we leave to future work.\\

\noindent \textsc{Remark 3.15} (dependence of $j_0$ on domain, wavelet system, $s$, and $\eta$) \phantomsection\label{Remark 3.15} Let $B\vcentcolon=\text{diam}(K)$. The proof of Lemma \hyperref[Lemma 3.14]{3.14} in fact shows that, for (\ref{eq3.13}) to hold, it is sufficient to choose $j_0\leqslant -F_s\log_2(G_{\psi,s}B/\eta)$ for some positive constants $F_s$ and $G_{\psi,s}$. In particular, $j_0$ scales logarithmically with $B/\eta$.\\

\noindent \textsc{Remark 3.16} (Simulations of $(W_{\text{wav}}^{(\psi,j_0)})_s$)\phantomsection\label{Remark 3.16} In Section \ref{subsec4.2}, we conduct the same simulations of translations and dilations as described in Remark \hyperref[Remark 3.10]{3.10} and Remark \hyperref[Remark 3.12]{3.12} but with our new wavelet $s$-Wasserstein distance formula. After avoiding the issues mentioned in the last subsection, the results are much better. In each fixed ``simple'' data set, after choosing a $j_0$ sufficiently small, $(W_{\text{wav}}^{(\psi,j_0)})_s/W_s$ is roughly a constant between all measures in the set.\\

\noindent \textsc{Remark 3.17} (Tightness of the lower and upper bounds in (\ref{eq3.13}))\phantomsection\label{Remark 3.17} When the domain containing the supports of the measures is large, a small $j_0$ is needed in the computation of $(W_{\text{wav}}^{(\psi,j_0)})_s$ to minimize the error $\eta$ in Theorem \hyperref[Theorem 3.13]{3.13}. However, with increasingly negative $j_0$, the difference between the lower and upper bound constants in (\ref{eq3.13}) becomes arbitrarily large. It is natural to ask if the bounds can actually be tightened. Unfortunately, we shall only give a partial answer --- that, the lower bound is indeed tight for any compact domain $K$ with a nonempty interior --- and leave the other half for future work. 

To see the tightness of the lower bound, fix $j_0$ and a wavelet system $(\phi,\psi)$. Since the domain has a nonempty interior, it contains an open ball $B_r$. Then, as $(\phi,\psi)$ is compactly supported, there exists $j'\in\mathbb{Z}$ and $\mathbf{k}'\in\mathbb{Z}^n$ so that $\psi_{(\epsilon,(j',\mathbf{k}'))}$ is supported in $B_r$ for all $\epsilon\in\alpha$. Now, pick $\epsilon'\in\alpha$ so that  $$||\psi_{(\epsilon',(j',\mathbf{k}'))}||_{(\dt{C}^s(\mathbb{R}^n))^{\ast}} = \max_{\epsilon\in\alpha}||\psi_{(\epsilon,(j',\mathbf{k}'))}||_{(\dt{C}^s(\mathbb{R}^n))^{\ast}}.$$
Since $\psi_{(\epsilon',(j',\mathbf{k}'))}$ integrates to $0$, the integrals of its positive and negative parts are equal, namely
$$\int_{\mathbb{R}^n}(\psi_{(\epsilon',(j',\mathbf{k}'))})^{+}d\mathbf{x} = \int_{\mathbb{R}^n}(\psi_{(\epsilon',(j',\mathbf{k}'))})^{-}d\mathbf{x} =C$$
for some constant $C>0$. Let $d\mu = ((\psi_{(\epsilon',(j',\mathbf{k}'))})^{+}/C)d\mathbf{x}$, $d\nu = ((\psi_{(\epsilon',(j',\mathbf{k}'))})^{-}/C)d\mathbf{x}\in\mathcal{P}_{\text{ab2}}(K)$, then we have, by definition,
\begin{align}
E_{\psi,s}W_s(\mu, \nu) &= \frac{1}{\max_{\epsilon\in\alpha}||\psi_{\epsilon}||_{(\dt{C}^s(\mathbb{R}^n))^{\ast}}}\sup_{f\in C^s_{1}(\mathbb{R}^n)}\int_{\mathbb{R}^n}f(\mathbf{x})\left(\frac{1}{C} \psi_{(\epsilon',(j',\mathbf{k}'))}(\mathbf{x})\right)d\mathbf{x}\notag\\
&= \frac{1}{C\max_{\epsilon\in\alpha}||\psi_{\epsilon}||_{(\dt{C}^s(\mathbb{R}^n))^{\ast}}}\max_{\epsilon\in\alpha}\sup_{f\in C^s_{1}(\mathbb{R}^n)}\Biggl|\int_{\mathbb{R}^n}f(\mathbf{x})\psi_{(\epsilon,(j',\mathbf{k}'))}(\mathbf{x})d\mathbf{x}\Biggr|\notag\\
&= \frac{2^{-j'(s+\frac{n}{2})}}{C\max_{\epsilon\in\alpha}||\psi_{\epsilon}||_{(\dt{C}^s(\mathbb{R}^n))^{\ast}}}\max_{\epsilon\in\alpha}||\psi_{\epsilon}||_{(\dt{C}^s(\mathbb{R}^n))^{\ast}}\notag\\
&= \frac{2^{-j'(s+\frac{n}{2})}}{C}\notag\\
&= \frac{2^{-j'(s+\frac{n}{2})}}{C}|(\psi_{(\epsilon',(j',\mathbf{k}'))})_{(\epsilon', (j', \mathbf{k}'))}^{(d)}|\notag\\
&=\sum_{\substack{\epsilon\in\alpha,\\ j\geqslant j_0,\  \mathbf{k}\in\mathbb{Z}^n}}2^{-j(s+\frac{n}{2})}\left|\left(\frac{\psi_{(\epsilon',(j',\mathbf{k}'))}}{C}\right)_{(\epsilon, (j, \mathbf{k}))}^{(d)}\right|\notag\\
&= (W_{\text{wav}}^{(\psi,j_0)})_s(\mu,\nu)\tag{eq3.17}\label{eq3.17}
\end{align}
In the second equality, we utilize our choice of $\epsilon'$ and the fact that $f$ and $-f$ exist in pairs in $C_1^s(\mathbb{R}^n)$; in the third equality, we follow the change of variable part of the proof in Lemma \hyperref[Lemma 3.4]{3.4}; and in the fifth and sixth equalities, we use the orthonormality of the wavelet system. Now, as the domain is fixed and the equalities in (\ref{eq3.17}) hold true for all $j_0\leqslant j'$, we may have arbitrarily small error $\eta$ as formulated in Theorem \hyperref[Theorem 3.13]{3.13}, and the lower bound constant has to be tight.

The upper bound constant is subtle and we leave the question open of its sharpness to future work. If one wants to find some particular $d\mu = p_{\mu}(\mathbf{x})d\mathbf{x}$, $d\nu = p_{\nu}(\mathbf{x})d\mathbf{x}\in\mathcal{P}_{\text{ab2}}(K)$ on certain $K$ such that, for some $j_0\in\mathbb{Z}$,
$$(W_{\text{wav}}^{(\psi,j_0)})_s(\mu,\nu) = \tilde{E}_{\psi,s,j_0}W_s(\mu,\nu),$$
then, one would need \begin{equation}\sum_{\substack{\epsilon\in\alpha,\\ j\geqslant j_0,\  \mathbf{k}\in\mathbb{Z}^n}}2^{-j(s+\frac{n}{2})}\text{sgn}((p_{\mu}-p_{\nu})_{(\epsilon, (j, \mathbf{k}))}^{(d)})\psi(\mathbf{x})_{(\epsilon, (j, \mathbf{k}))}\notag\end{equation}
to precisely be a $f_{\text{max}}$ that achieves the supremum of 
$$\sup_{f\in \dt{C}^s_{\tilde{E}_{\psi,s,j_0}}(K)}\int_{K}f(\mathbf{x}) (p_{\mu}(\mathbf{x})-p_{\nu}(\mathbf{x}))d\mathbf{x}.$$
where, particularly for negative $j_0$, 
\begin{align}
\tilde{E}_{\psi,s,j_0} = &\frac{\sum_{i=1}^n||\textstyle \sum\psi^{(i)}||_{\infty}}{2^{1-s}(2^{1-s}-1)}+\frac{2||\textstyle\sum\psi||_{\infty}}{1-2^{-s}} + \notag\\
&|j_0|\sup_{\substack{\mathbf{x}\neq\mathbf{y}\\\eta_k\in\{+1,-1\}}}\left|\frac{\sum_{\mathbf{k}\in\mathbb{Z}^n,\epsilon\in\alpha}\eta_k\psi_{(\epsilon,(0,\mathbf{k}))}(\mathbf{x})-\sum_{\mathbf{k}\in\mathbb{Z}^n,\epsilon\in\alpha}\eta_k\psi_{(\epsilon,(0,\mathbf{k}))}(\mathbf{y})}{|\mathbf{x}-\mathbf{y}|^s} \right|\notag
\end{align}
is much harder to handle than the lower bound constant. It is unclear how one can in general construct an example for the upper bound.

Numerically speaking, it does not appear that the lower and upper bounds comparing the new wavelet and classical $s$-Wasserstein distances can be tightened to be roughly equal. In Section \ref{subsec4.2}, the domains of measures in each of our data set is only moderately large and our choices of $j_0$ are roughly capable of making the $\eta$ errors, which appear in Theorem \hyperref[Theorem 3.13]{3.13}, comparably small. We have already proved the lower bound constant $E_{\psi,s}$ is tight and neither depends on the domain $K$ nor $j_0$. If the upper bound constant $\tilde{E}_{\psi,s,j_0}$ were really close to the lower one $E_{\psi,s}$, then, after fixing the wavelet system $(\phi,\psi)$ and $s$, in principle, we could choose almost the same normalizing constant for all of our data sets when comparing to the classical Wassersetin $s$-distance. This was not quite the case numerically. Furthermore, numerical simulations suggest that the upper bound constant increases with $j_0$ at least, as in (\ref{eq3.13}); see Simulation \hyperref[Simulation seven]{seven} where we consider the behavior as $j_0$ varies.

Of course, one might still ask, just for translations and dilations of a fixed type of common measure, under fixed $j_0$, whether we can expect some tightened bounds like we numerically observe in Section \ref{subsec4.2}, which is essentially the reason why we give a high evaluation to the numerical performance of our new wavelet $s$-Wasserstein distance. This direction is also left for future work. \\

\noindent \textsc{Remark 3.18} (Extension of Theorem \hyperref[Theorem 3.13]{3.13} to a larger family of probability measures)\phantomsection\label{Remark 3.18}
In this remark, we discuss some of our observations regarding whether the equivalence result in Theorem \hyperref[Theorem 3.13]{3.13} be extended to a larger family of probability measures than $\mathcal{P}_{\text{ab2}}(K)$. In below, we naturally extend the definition (\ref{eq2.8}) to probability measures $\mu$ and use $\mu_{(\epsilon, (j, \mathbf{k}))}^{(d)}$ to denote the detail coefficients.

One potential approach to extending our result could proceed as follows: first, $\mathcal{P}_{\text{ab2}}(K)$ is dense in $\mathcal{P}_{s}(K)$ in the $s$-Wasserstein distance on any compact $K$ by \cite[Remark 5.13]{bib3} and Jensen's inequality; second, the Wasserstein distance topology on $\mathcal{P}_{s}(K)$ is precisely the weak-$\ast$ topology (against compactly supported continuous function) when the domain is compact \cite[Theorem 5.10]{bib3}; third, the constants in our equivalence (\ref{eq3.13}) do not depend on a priori estimate on the densities of measures. Combining them, it is tempting to extend the claim in Theorem \hyperref[Theorem 3.13]{3.13} to all $\mathcal{P}_{s}(K)$ directly through a density argument. However, we ultimately have to face the problem of exchanging the limit operation in the approximation and the infinite sums in defining $(W_{\text{wav}}^{(\psi,j_0)})_s$, where a dominating sequence in $\ell^1$ is required to fully justify this procedure. Without some additional requirement on the wavelet coefficients of the measures (like the existing $\ell^2$ constraint), it is hard to complete the proof.

For a second potential approach, given $\mu,\nu\in\mathcal{P}_s(K)$ and a wavelet system $(\phi,\psi)$, one could try to work with the measures arising from the partial sums $$\mu_{j_0}^{j'}\vcentcolon = \sum_{\substack{\epsilon\in\alpha,\\ j_0\leqslant j\leqslant j',\  \mathbf{k}\in\mathbb{Z}^n}}2^{-j(s+\frac{n}{2})}\mu_{(\epsilon, (j, \mathbf{k}))}^{(d)}\psi(\mathbf{x})_{(\epsilon, (j, \mathbf{k}))}d\mathbf{x}$$
and 
$$\nu_{j_0}^{j'}\vcentcolon = \sum_{\substack{\epsilon\in\alpha,\\ j_0\leqslant j\leqslant j',\  \mathbf{k}\in\mathbb{Z}^n}}2^{-j(s+\frac{n}{2})}\nu_{(\epsilon, (j, \mathbf{k}))}^{(d)}\psi(\mathbf{x})_{(\epsilon, (j, \mathbf{k}))}d\mathbf{x},$$
with $j_0,j'\in\mathbb{Z}$ and $j_0\leqslant j'$. By the compactly supported nature of the measures and the wavelet system, $\mu_{j_0}^{j'}$ and $\nu_{j_0}^{j'}$ consist of finitely many terms. Their difference may not be the difference between two probability measures, but has zero total mass and a square-integrable density. After scaling, our existing equivalence result (\ref{eq3.13}) implies that, given an error $\eta$, 
\begin{equation}
E_{\psi,s}||\mu_{j_0}^{j'}-\nu_{j_0}^{j'}||_{(\dt{C}^s(\mathbb{R}^n))^{\ast}}-\eta\leqslant (W_{\text{wav}}^{(\psi,j_0)})_s(\mu_{j_0}^{j'}, \nu_{j_0}^{j'}) \leqslant \tilde{E}_{\psi,s,j_0}||\mu_{j_0}^{j'}-\nu_{j_0}^{j'}||_{(\dt{C}^s(\mathbb{R}^n))^{\ast}}\notag
\end{equation}
holds for all $j'\in\mathbb{Z}$ when $j_0\in\mathbb{Z}$ is sufficiently small.
We then proceed with a similar argument as how we prove the lower bound in Theorem \hyperref[Theorem 3.13]{3.13}: leave an epsilon of room to connect $||\mu_{j_0}^{j'}-\nu_{j_0}^{j'}||_{(\dt{C}^s(\mathbb{R}^n))^{\ast}}$ and $W_s(\mu, \nu)$; this obviously requires a slightly enhanced version of Lemma \hyperref[Lemma 3.14]{3.14}, but such modified Lemma can be proved. As shown in (\ref{eq3.15}), the fundamental issue of this type of argument is that a uniform $L^1$ type upper bound is needed for $\mu_{j_0}^{j'}-\nu_{j_0}^{j'}$, for all $j'$ large. Unfortunately, one can check that such bound is not achievable if $\mu$ is the Dirac delta measure supported at the origin. Using this method, however, it is still possible to extend the equivalence results to allow measures with densities in the (real) Hardy space $\mathcal{H}^1(\mathbb{R}^n)$; see \cite[Chapter 5.2]{bib4} for a definition.

We have intentionally been sketchy here for two reasons: first, there might be other more involved arguments to extend the equivalence maximally to $\mathcal{P}_s(K)$. Second, even for the methods we tried, there could exist other spaces that the arguments work. Practically, the existing framework on $\mathcal{P}_{\text{ab2}}(K)$ is sufficient for this paper and the ultimate determination of the problem of extension shall be left for future works. \\

\noindent \textsc{Remark 3.19} (Dependence of $(W_{\text{wav}}^{(\psi,j_0)})_s$ on $j_0$)\phantomsection\label{Remark 3.19} Another natural question to ask with our new wavelet $s$-Wasserstein distance is whether there is a pattern of rate of increase in $j_0$ when computing the distance between two arbitrary probability measures in $\mathcal{P}_{\text{ab2}}(K)$ for some $K$. The answer is a clear no, as illustrated in the following example, in which measures $\mu$ and $\nu$ are constructed in the same spirit as in remark \hyperref[Remark 3.17]{3.17}, so that $(W_{wav}^{\psi,j_0})_s(\mu,\nu)$ increases arbitrarily with $j_0$.

Fix a system $(\phi,\psi)$, select any $j'\in\mathbb{Z}$, $j_0\leqslant j'$, $\epsilon'\in\alpha$, and pick arbitrary nonnegative function $a(j)$ on $\mathbb{Z}$. Consider
$$g(\mathbf{x}) = \sum_{j=j_0}^{j'}a(j)2^{j(s+n)}\psi_{\epsilon'}(2^j\mathbf{x})$$
Since $\psi_{\epsilon'}$ integrates to 0 and is compacted supported, we can have $d\mu = (g(\mathbf{x})^{+}/C)d\mathbf{x}$, $d\nu = (g(\mathbf{x})^{-}/C)d\mathbf{x}\in\mathcal{P}_{\text{ab2}}(K)$ for some constant $C$ on some $K$. Then, we have
\begin{align}(W_{\text{wav}}^{(\psi,j_0)})_s(\mu,\nu) &= \frac{1}{C}\sum_{\substack{\epsilon\in\alpha,\\ j\geqslant j_0,\  \mathbf{k}\in\mathbb{Z}^n}}2^{-j(s+\frac{n}{2})}|(g(\mathbf{x}))_{(\epsilon, (j, \mathbf{k}))}^{(d)}|\notag\\
&= \frac{1}{C}\sum_{\substack{\epsilon\in\alpha,\\ j\geqslant j_0,\  \mathbf{k}\in\mathbb{Z}^n}}2^{-j(s+\frac{n}{2})}\Biggl|\Biggl(\sum_{l=j_0}^{j'}a(l)2^{l(s+n)}\psi_{\epsilon'}(2^l\mathbf{x})\Biggr)_{(\epsilon, (j, \mathbf{k}))}^{(d)}\Biggr|\notag\\
&= \frac{1}{C}\sum_{j=j_0}^{j'}2^{-j(s+\frac{n}{2})}a(j)2^{j(s+\frac{n}{2})}\notag\\
&= \frac{1}{C}\sum_{j=j_0}^{j'}a(j)\notag
\end{align}
by the orthonormality of the wavelet system. As $a(j)$ is an arbitrary nonnegative function, there is no uniform pattern on how $(W_{\text{wav}}^{(\psi,j_0)})_s$ is going to increase with respect to the decrease of $j_0$ for different measures in $\mathcal{P}_{\text{ab2}}(K)$. 

Note however that this example does not immediately tell us how the bounds in (\ref{eq3.13}) should behave, since it is hard to estimate $W_s(\mu,\nu)$. 

\subsection{Wavelet linearized optimal transport}\label{subsec3.5}
In this subsection, we define a map with a goal similar to the one discussed by Wang et al. in \cite{bib10}. This work introduced an embedding of $\mathcal{P}_{2}(\mathbb{R}^n)$ into a linear space, known as the \textbf{linearized optimal transport (LOT)}, which is inspired by the formal Riemannian structure of $W_2$ and preserves certain key properties. 

We begin by recalling the LOT embedding. In the LOT framework, one first fixes a reference measure $\sigma\in\mathcal{P}_{2}(\mathbb{R}^n)$ that is absolutely continuous to the Lebesgue measure. Then, by   Brenier's theorem \cite[Theorem 1.22]{bib3}, for any $\mu\in\mathcal{P}_{2}(\mathbb{R}^n)$, there exists an \textbf{optimal transport map} $T_{\sigma}^{\mu}: \mathbb{R}^n\to\mathbb{R}^n$ such that 
\begin{equation}
W_2(\mu, \nu) = \Biggl(\int_{\mathbb{R}^n}|T_{\sigma}^{\mu}(\mathbf{x})-\mathbf{x}|^2d\sigma(\mathbf{x})\Biggr)^{\frac{1}{2}}.\notag
\end{equation}
The LOT embedding map is then defined as $\mathcal{T}_{\sigma}: \mathcal{P}_{2}(\mathbb{R}^n) \to L^2_{\sigma}(\mathbb{R}^n)$  
$$\mathcal{T}_{\sigma}: \mu\mapsto T_{\sigma}^{\mu},$$
where $L^2_{\sigma}(\mathbb{R}^n)$ is the space of square-integrable functions against $\sigma$.

The map $\mathcal{T}_{\sigma}$ is by definition injective; also, since $(T_{\sigma}^{\mu}, T_{\sigma}^{\nu})_{\#}\sigma$ is a transport plan between $\mu$ and $\nu$, by (\ref{eq2.1})
\begin{equation}W_2(\mu,\nu)\leqslant ||T_{\sigma}^{\mu} - T_{\sigma}^{\nu}||_{L_{\sigma}^2(\mathbb{R}^n)}.\notag
\end{equation}
The right-hand side of the above equation is defined to be the linearized 2-Wasserstein distance $(W_{\text{LOT},\sigma})_2(\mu,\nu)$ between $\mu$ and $\nu$. Under appropriate hypotheses on the reference measure $\sigma$, the underlying domain, and $\mu$, $\nu$, Delalande and M{\'e}rigot showed that $(W_{\text{LOT},\sigma})_2$ is bi-Hölder equivalent to $W_2$ (Delalande \& M{\'e}rigot \cite[Corollary 3.4, Theorems 4.2, and 4.3]{bib11}). Importantly, $W_2(\mu,\nu)$ and $(W_{\text{LOT},\sigma})_2(\mu,\nu)$ coincide when $\mu$ and $\nu$ are just translations and dilations of each other (\cite[Moosmüller \& Cloninger]{bib12}) and only produce small errors when there are perturbations.

The main motivation for introducing such a linear embedding is the following: previously, in a data set of $N$ measures, to obtain all pairwise 2-Wasserstein distances, $N(N-1)/2$ number of optimal transport maps are needed, which are all expensive to compute. With linearized optimal transport, only $N$ optimal transport maps between the reference measure and the measures in the data set need to be computed, and then $N(N-1)/2$ cheap $L^2$ norm computations yield all pairwise linearized 2-Wasserstein distances. Moreover, in the embedded linear space, one may apply many supervised and unsupervised learning methods not available in nonlinear spaces.

We now cast our result into a wavelet alternative to the framework above. Based on the development in the previous subsections, we only formulate the embedding with our new wavelet $s$-Wasserstein distance for $0<s<1$ in Theorem \hyperref[Theorem 3.13]{3.13}, though an analogous approach can be used for the original wavelet $s$-Wasserstein distance.

After fixing a compact domain $K$ and then $j_0$ sufficiently small accordingly, define the map: $\mathcal{T}_{\text{wav}}^{j_0}: \mathcal{P}_{\text{ab2}}(K)\to \ell^2\left(\alpha\bigoplus \{j\geqslant j_0\}\bigoplus\mathbb{Z}^n\right)$ by
$$\mathcal{T}_{\text{wav}}^{j_0}: p_{\mu}(\mathbf{x})d\mathbf{x}\mapsto((p_{\mu})_{(\epsilon, (j,\mathbf{k}))}^{(d)})_{\substack{\epsilon\in\alpha,\\j\geqslant j_0,\mathbf{k}\in\mathbb{Z}^n}}.$$
If we equip $\mathcal{P}_{\text{ab2}}(K)$ with the classical $W_s$ and equip the codomain $\ell^2(\alpha\bigoplus \{j\geqslant j_0\}$$\bigoplus\mathbb{Z}^n)$ of the map $\mathcal{T}_{\text{wav}}^{j_0}$ with the weighted $\ell^1$ norm $$\Bigl|\Bigl|(p_{(\epsilon,(j,\mathbf{k}))}^{(d)})_{\substack{\epsilon\in\alpha,\\j\geqslant j_0,\mathbf{k}\in\mathbb{Z}^n}}\Bigr|\Bigr|\vcentcolon = \sum_{\substack{\epsilon\in\alpha,\\ j\geqslant j_0,\  \mathbf{k}\in\mathbb{Z}^n}}2^{-j(s+\frac{n}{2})}|p_{(\epsilon, (j, \mathbf{k}))}^{(d)}|,$$
then, similar to LOT, $\mathcal{T}_{\text{wav}}^{j_0}$ is an embedding map to a linear space with an equivalent distance (up to an error) as shown in Theorem \hyperref[Theorem 3.13]{3.13}. We will call this embedding the \textbf{wavelet linearized optimal transport (WLOT)}.

There are some immediate differences between the classical LOT and WLOT, aside from the parameter $s$ in the $s$-Wasserstein distance. Clearly, the latter does not require a choice of reference measure and directly embeds the target measures instead of the optimal transport maps between the reference measure and target measures. The transport map might be complicated and losses certain important features of the original target measures, like the sparsity, which are important in certain situations like image classification. In contrast, those properties are largely preserved by our wavelet coefficient embedding, due to the localization properties of the wavelet systems.

We close by comparing our approach to another linear embedding for the $2$-Wasserstein distance, which utilizes the dual norm of a homogenous Sobolev space and, similar to ours, directly embed the target measures. Let \(\dot{H}^1(\mathbb{R}^n)\) be the quotient of the completion of smooth, compactly supported functions \(C_c^\infty(\mathbb{R}^n)\) under the seminorm$$\| f \|_{\dot{H}^1(\mathbb{R}^n)} = \| \nabla f \|_{L^2(\mathbb{R}^n)},$$
such that it becomes a complete normed space. Then, the space \(\dot{H}^{-1}(\mathbb{R}^n)\) can be viewed as its topological dual, which includes as elements the difference between any pair of probability measures on $\mathbb{R}^n$, which has zero total mass. Peyre proved that $W_2$ is equivalent to $\dot{H}^{-1}$ on subsets of probability measures on $\mathbb{R}^n$ that possess finite second moment and densities that are uniformly bounded from below and above, with the constants in the equivalence depending on those bounds and blowing up at extremes \cite{bib27}. The dependence of the constants on the attributes of the target measures themselves makes such embedding ineffective in dealing with scattered measures with low densities or concentrated measures with high densities, which are common in practical data sets. In contrast, the constants in the equivalence of our embedding only depend on the domain of target measures.

In Greengard et al. \cite{bib28}, it was also proved that certain weighted $\dot{H}^{-1}$ norms provide a local linearization of $W_2$, with quantitative estimates depending on the measure around which one is linearizing. However, it is a local result around a measure and, in essence, for measures differ only by perturbations. Our theorem is global.

\section{Numerical simulations}\label{sec4}
\subsection{Numerically computing the wavelet $s$-Wasserstein distances}\label{subsec4.1}
\noindent We now explain the details of our numerical approach for computing the new wavelet $s$-Wasserstein distance $(W_{\text{wav}}^{(\psi,j_0)})_s$ for $0<s\leqslant1$. The other two formulations appearing in this paper can be computed analogously; see Remark \hyperref[Remark 4.3]{4.3}. In spite of the theoretical limitations in the $s=1$ case described in the previous section, we include simulations with $s=1$ in view of the natural continuity for the classical $W_s$ distance as $s\to1$ in Proposition \hyperref[Proposition 2.5]{2.5}. Our simulations are implemented using \texttt{Python}, with the help of the \texttt{numpy} \cite{bib14} and \texttt{pywt} (PyWavelets) \cite{bib15} packages. As our main goal is to explore the qualitative properties of various wavelet $s$-Wasserstein distances and compare their behavior to the classical $s$-Wasserstein distance (Section \ref{subsec4.2}), we focus on probability measures on the real line $\mathbb{R}$.

In each simulation, we begin with a finite set of functions $\{p_{i}\}_{i\in \llbracket m\rrbracket}$ and a compact set $K \subset \mathbb{R}$ such that $p_{i}(x)dx\in\mathcal{P}_{\text{ab2}}(K)$  for all $i\in \llbracket m\rrbracket$. \\ 

\noindent \textbf{Step one: preparation}. Use the \texttt{pywt.Wavelet} function to pick a wavelet system that satisfies the requirement in Theorem \hyperref[Theorem 3.13]{3.13}. Select $j_0\in\mathbb{Z}$ small enough such that it reduces the error $\eta$ as in Theorem \hyperref[Theorem 3.13]{3.13} and so that roughly $B:={\rm diam}(K)\leq 2^{-j_0}$. In view of Remark \hyperref[Remark 2.3]{2.3}, one may translate functions $\{p_{i}\}_{i\in \llbracket m\rrbracket}$ to be supported on some the interval $[0,B]$. Then, assume without loss of generality that ${\rm supp}(\psi) \subseteq [0,+\infty)$, which is satisfied by almost all constructed wavelet systems. All of these are to ensure that, for all $j\geqslant j_0$, the coefficients  $(p_i)_{(\epsilon,(j,\mathbf{k}))}^{(d)}$ vanish whenever   $k<0$ or  $k>2^{j-j_0}$. \\

\noindent \textbf{Step two: initialization}. Fix $M\in\mathbb{N}$ sufficiently large (with $M>|j_0|$ if $j_0<0$). For any two functions, say $p_1$ and $p_2$, in the collection, initialize $p_1-p_2$ as in Section \ref{subsec2.4}, (\ref{eq2.13}); it amounts to sampling $p_1-p_2$ on a uniform grid on the interval $[0,2^{-j_0}]$, with grid spacing $2^{-(j_0+M)}$, and multiplying the values by a factor $2^{-(j_0+M)/2}$. This yields an initialization vector (list) of length $2^M$, which we denote by $\mathbf{A}^{j_0,M}$. In this way, the size of $M$ controls the accuracy of initialization, as well as summation in Step four below. \\

\noindent \textbf{Step three: discrete wavelet transformation}. Input the initialization vector $\mathbf{A}^{j_0,M}$ above into the \texttt{pywt.wavedec} function, together with the chosen wavelet system, the \texttt{"zero"} extension mode\footnote{It is the most natural mode in our case, as the measures should simply vanish outside of their respective supports.}, and \texttt{level} = $M$.\\

\noindent \textbf{Step four: summation}. In \texttt{pywt}, conveniently, the output of the previous step groups the coefficients in a level-wise fashion: instead of being precisely like (\ref{eq2.12}), the result is (with parentheses denoting the starting and ending of arrays)\footnote{A few more coefficients will be computed at each level, which we do not explicitly write out for convenience; see Remark \hyperref[Remark 2.7]{2.7}. It does not affect the logic flow of the method.}
\begin{align}
\mathbf{O}^{j_0,M} = \bigl(\bigl((p_1-p_2)_{(j_0,0)}^{(a)}\bigr),&\bigl((p_1-p_2)_{(j_0,0)}^{(d)}\bigr), \bigl((p_1-p_2)_{(j_0+1,0)}^{(d)}, (p_1-p_2)_{(j_0+1,1)}^{(d)}\bigr),\bigl(\cdots\bigr),\notag\\
&\bigl((p_1-p_2)_{(j_0+M-1,0)}^{(d)},\cdots,(p_1-p_2)_{(j_0+M-1,2^{M-1}-1)}^{(d)}\bigr)\bigr).\notag
\end{align}
The structure of this output allows one to easily apply by the weight $2^{-j(s+1/2)}$ at each level, $j=j_0,\cdots,j_0+M-1$. Then, discard the first approximation coefficient and sum the absolute values of all detail coefficients with weights. Due to the preparation in step one, the final value is
\begin{equation}\sum_{j= j_0}^{j_0+M-1}\sum_{k=0}^{2^{j-j_0}}2^{-j(s+\frac{1}{2})}|(p_{1}-p_{2})_{(j, k)}^{(d)}| = \sum_{j= j_0}^{j_0+M-1}\sum_{k\in\mathbb{Z}}2^{-j(s+\frac{1}{2})}|(p_{1}-p_{2})_{(j, k)}^{(d)}|.\tag{eq4.1}\label{eq4.1}\end{equation}
Since the sum in $(W_{\text{wav}}^{(\psi,j_0)})_s(p_1dx,p_2dx)$ converges, the numerical result (\ref{eq4.1}), which we denote as $(W_{\text{wav,app}}^{(\psi,j_0, M)})_s(p_1dx,p_2dx)$, provides a good approximation to $(W_{\text{wav}}^{(\psi,j_0)})_s(p_1dx,p_2dx)$ for $M$ large. \\

\noindent \textbf{Step five: repetition}: Repeat from step two for all pairs of functions in the set $\{p_{i}\}_{i\in \llbracket m\rrbracket}$ whose distances are wished to be computed. At last, when all the desired values are obtained, if one wants to make the result comparable in magnitude with the classical $W_s$, multiply them by a uniform normalizing factor $C_{(\psi,s,j_0,M)}^{\text{nor},\{p_{i}\}_{i\in \llbracket m\rrbracket}}$.\\

\noindent \textsc{Remark 4.1} (Order of complexity)
With a detailed recount of the procedure, it is clear that this computation is linear in time with respect to the number of samples $2^M$ above collected on a given difference function $p_1-p_2$: the wavelet transform \texttt{pywt.wavedec} itself and the summation of the weighted coefficients are both linear in time, with respect to the inputs, while the time complexity of the other operations in the steps, namely multiplication by a uniform factor $2^{-(j_0+M)/2}$ in initialization and application of weights $2^{-j(s+1/2)}$ at different levels in the final summation, will be majorized by the previous two when $M$ is large.\\

\noindent \textsc{Remark 4.2} (The normalizating factor $C_{(\psi,s,j_0,M)}^{\text{nor},\{p_{i}\}_{i\in \llbracket m\rrbracket}}$) \phantomsection\label{Remark 4.2} In practice, the constants in the upper and lower bounds in Theorem \hyperref[Theorem 3.13]{3.13}, in which we compare $(W_{\text{wav}}^{(\psi, j_0)})_s$ to $W_s$, may be far away from 1. Consequently, to visually analyze the behavior of the two distances on the same plot, we multiply $(W_{\text{wav}}^{(\psi, j_0)})_s$  by a normalizing constant, which depends on the collection of functions $\{ p_i \}_{i \in \llbracket m\rrbracket}$ that we consider. As our theoretical derivation does not give guidance on the selection of this factor, we select it empirically, based on the observations of the values of $(W_{\text{wav,app}}^{(\psi,j_0, M)})_s(p_idx,p_jdx)$ for all $p_i, p_j$ in the collection whose distances need to be computed.  We mention that, in many applications, only relative magnitude of the distance matters. In this case, one does not even need to worry about picking this factor.\\ 

\noindent \textsc{Remark 4.3} (Computations for other wavelet $s$-Wasserstein distances)\phantomsection\label{Remark 4.3}
In the next subsection, we shall also display some numerical approximation results, $(W_{\text{wav,ori,app}}^{(\psi,j_0,M)})_s$ and $(W_{\text{wav,ori,app}}^{(\phi,\psi,C_0,C_1,j_0,M)})_s$, for the two formulations of the original wavelet $s$-Wasserstein distances (the former being the $C_0 = 0$ and $C_1 = 1$ case of the later). The steps for computing them are just slightly different from above. First, the discrete wavelet transform in step three should now be decomposed with $\texttt{level}=j_0+M$, but not fully to $\texttt{level}=M$, to obtain the approximation coefficients at $j=0$; second, the final summation formula in step four needs to be adjusted according to Theorem \hyperref[Theorem 3.7]{3.7} and Theorem \hyperref[Theorem 3.11]{3.11}. Needless to say, the normalzation constants $C_{(\psi,s,j_0,M)}^{\text{ori,nor},\{p_{i}\}_{i\in \llbracket m\rrbracket}}$, $C_{(\phi,\psi,C_0,C_1,s,j_0,M)}^{\text{ori,nor},\{p_{i}\}_{i\in \llbracket m\rrbracket}}$ will, in general, be different. We remark lastly that, although choosing $j_0$ and $M$ are not theoretically required for those two distances, it is still necessary from a computational perspective when using DWT. This is because only with a sufficiently small $j_0$ and a positive, large enough $M$, the DWT output will give all nonzero coefficients used in those two distances for levels below $j_0+M$.

\subsection{Numerical comparison between wavelet $s$-Wasserstein distances and $s$-Wasserstein distance}\label{subsec4.2}
We display below some numerical simulations of the three formulations of the wavelet $s$-Wasserstein distances for $0<s\leqslant1$ and compare them to $W_s$ when the measures are translations or dilations of each other. 

In this subsection, the wavelet system $(\phi,\psi)$ and the domain $K$ are chosen to be the Daubechies 10 wavelet (\texttt{db10}) (\cite[Chapter 6,7]{bib7}) and $[0,3]$, respectively. The $s$ values are $1,1/2, 1/4$. For the alternative formulation of the original wavelet $s$-Wasserstein distance in Theorem \hyperref[Theorem 3.11]{3.11}, the constant $C_0$ is selected to be $3^s$ and $C_1$ to be 1. To compare with the classical $s$-Wasserstein distance $W_s$, we use the OT Network Simplex solver \texttt{ot.emd} in the \texttt{pot} library \cite{bib16}, input with the discretizations of the data functions obtained by sampling them uniformly with 1000 points on the domain.\\

\noindent \textsc{Remark 4.4} (Notational conventions on plots) Before starting, we fix some notational conventions used particularly on the plots displayed in both this and next subsections. First, in our first reformulation, which appeared in Theorem \hyperref[Theorem 3.11]{3.11}, we technically still use the formula for original $s$-Wasserstein distance but enforce the choice of $C_0$ to be greater 0. However, to distinguish this reformulation from the original one with the particular choice $C_0 = 0$ and $C_1=1$, we change the subscript of this distance in the reformulation from ``ori'' to ``alt''.

Next, due to the abundance of super(sub)scripts in the normalizing constants for different wavelet distances, we only label the variables that change on each plot and suppress all the others. Similar conventions are also used for the wavelet distances themselves. However, ``wav",``ori", and ``alt" are not ignored but abbreviated as ``w", ``o", and ``a" respectively. Therefore, for example, if we are displaying the first reformulation of the original wavelet $s$-Wasserstein distance, which appeared in Theorem \hyperref[Theorem 3.11]{3.11}, on a plot where only the parameter $s$ is varied, then this quantity will be labeled as $C_s^{\text{a}}(W_{\text{w,a}})_s$ instead of $C_{(\phi,\psi,C_0,C_1s,j_0,M)}^{\text{alt,nor},\{p_{i}\}_{i\in \llbracket m\rrbracket}}(W_{\text{wav,alt,app}}^{(\phi,\psi,C_0,C_1,j_0,M)})_s$. One should remember, however, all the dependencies on various quantities, and be careful when comparing results across different plots. Lastly, the result from $\texttt{pot}$ will be directly labeled as $W_s$.\\
\subsubsection{Translations and dilations of the uniform probability measure}
\noindent \textbf{Simulation one}. 
\phantomsection\label{Simulation one}
(Translations of the uniform probability measure on $[0,1]$) Let $p(x) \vcentcolon= \mathbbm{1}_{[0,1]}(x)$ be the indicator function supported $[0,1]$ and $p_{a}(x)\vcentcolon = p(x-a)$ be the translation of $p$ by $a>0$. Three wavelet $s$-Wasserstein distances between $p$ and $p_{a}$ for 20 evenly spaced values of  $a\in[0,2]$ are computed for three different choices of $s$, with comparison to the classical $s$-Wasserstein distance. The $j_0$ and $M$ are chosen to be $-11$ and $22$.

\begin{figure}[H]
\phantomsection\label{Fig1(a)}
\centering
\subfloat[\centering $C_{s}^{\text{o}}(W_{\text{w,o}})_s$ with $C_{s}^{\text{o}}  = 1/4, 1/5, 1/6$ for $s=1, 0.5, 0.25$ respectively and $W_s$]{{\includegraphics[width=10cm]{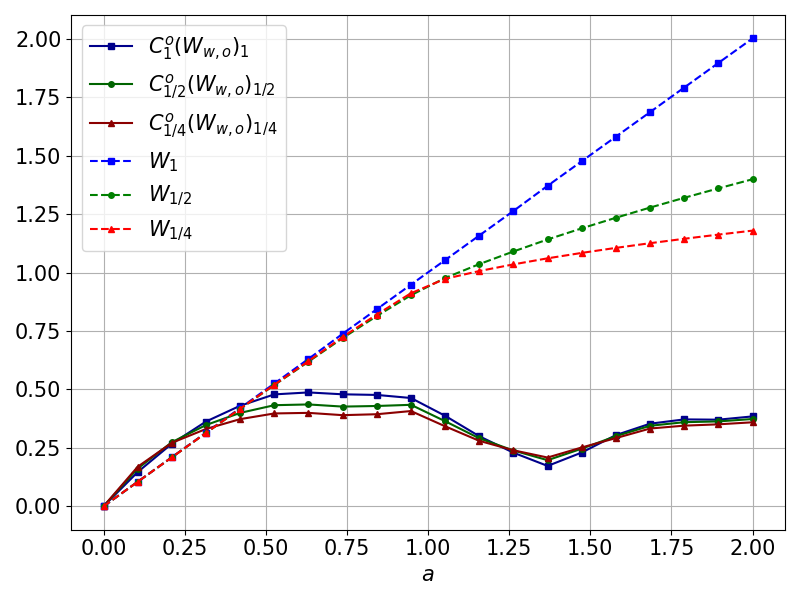}}}%
\end{figure}
\begin{figure}[H]
\phantomsection\label{Fig1(b)}
\centering
\subfloat[\centering $C_{s}^{\text{a}}(W_{\text{w,a}})_s$ with $C_{s}^{\text{a}}  = 1/8, 2/13, 1/7$ for $s=1, 0.5, 0.25$ respectively and $W_s$]{{\includegraphics[width=10cm]{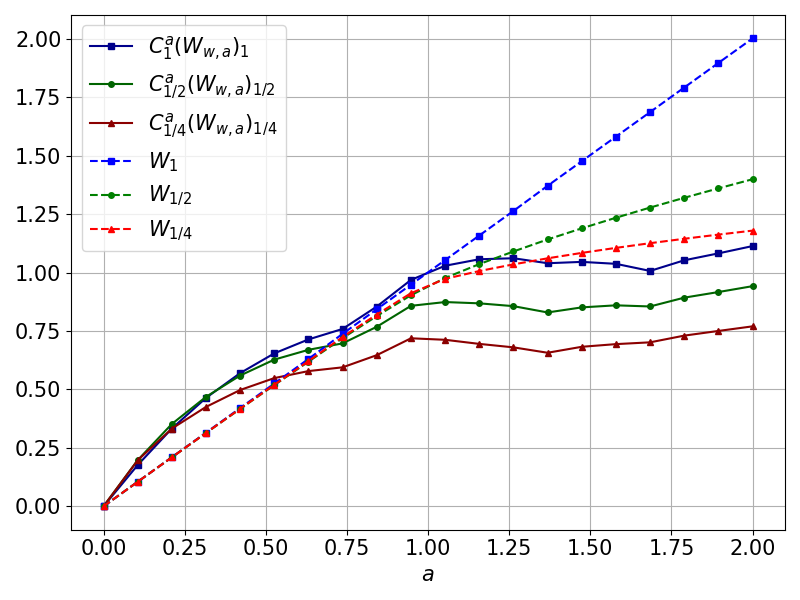} }}%
\end{figure}
\begin{figure}[H]
\phantomsection\label{Fig1(c)}
\centering
\subfloat[\centering $C_{s}(W_{\text{w}})_s$ with $C_{s} = 1/110, 1/22, 4/55$ for $s=1, 0.5, 0.25$ respectively and $W_s$]{{\includegraphics[width=10cm]{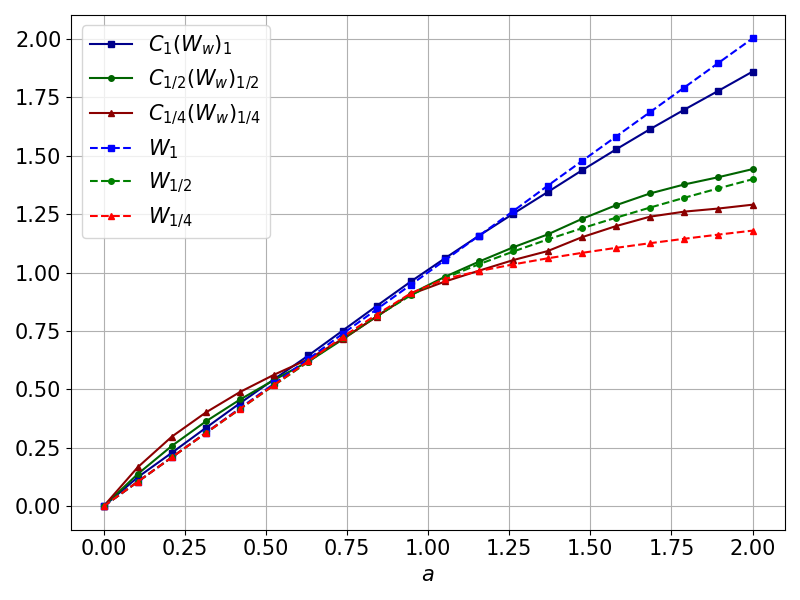} }}%
    \caption{Normalized numerical results of three wavelet $s$-Wasserstein distances in comparison with the classical $s$-Wasserstein distance for translations of the uniform probability measure on $[0,1]$}
    \label{Fig1}
\end{figure}
Both \hyperref[Fig1(a)]{Fig. 1 (a)} and \hyperref[Fig1(b)]{Fig. 1 (b)} show the stagnating, nonincreasing behavior of the original $s$-Wasserstein distance (one original and the other one reformulated) when the translations between uniform probability measures are large, demonstrating the incompatibility with translations, no matter what normalization is chosen; see Remark \hyperref[Remark 3.10]{3.10} and Remark \hyperref[Remark 3.12]{3.12}. \hyperref[Fig1(c)]{Fig. 1 (c)} provides the best performance: even though our main Theorem \hyperref[Theorem 3.13]{3.13} shows that the new wavelet $s$-Wasserstein distance is merely equivalent to the classical $s$-Wasserstein distance, up to an error, we see numerically that the behavior of $W_s$ under translations is well captured.\\

\noindent \textbf{Simulation two}. (Dilations of the uniform probability measure on $[1,2]$) Let $p(x) \vcentcolon= \mathbbm{1}_{[1,2]}(x)$ be the indicator function supported $[1,2]$ and $p_{b}(x)\vcentcolon = \mathbbm{1}_{[-(1/2)b+3/2,(1/2)b+3/2]}(x)$ be the dilation of $p$ around its center $3/2$ by $b>0$. Three wavelets $s$-Wasserstein distances between $p$ and $p_{b}$ for  20 evenly spaces values of $b\in[1/2,3/2]$ are computed for three different $s$ and compared to the classical $s$-Wasserstein distance. The $j_0$ and $M$ are chosen to be $-9$ and $22$.\\
\begin{figure}[H]
\phantomsection\label{Fig2(a)}
\centering
\subfloat[\centering $C_{s}^{\text{o}}(W_{\text{w,o}})_s$ with $C_{s}^{\text{o}} = 2/17, 2/17, 2/17$ for $s=1, 0.5, 0.25$ respectively and $W_s$]{{\includegraphics[width=10cm]{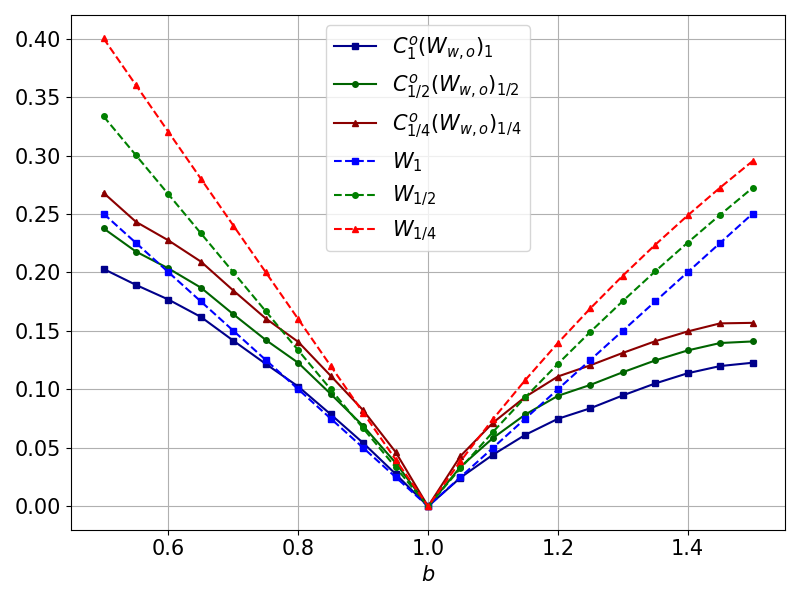} }}%
\end{figure}
\begin{figure}[H]
\phantomsection\label{Fig2(b)}
\centering
\subfloat[\centering $C_{s}^{\text{a}}(W_{\text{w,a}})_s$ with $C_{s}^{\text{a}} = 1/16, 1/11, 1/10$ for $s=1, 0.5, 0.25$ respectively and $W_s$]{{\includegraphics[width=10cm]{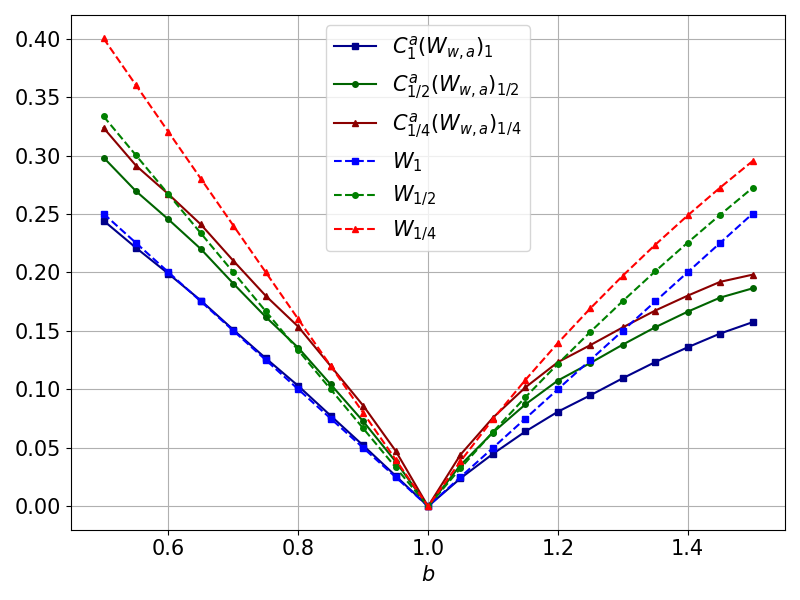} }}%
\end{figure}
\begin{figure}[H]
\phantomsection\label{Fig2(c)}
\centering
\subfloat[\centering $C_{s}(W_{\text{w}})_s$ with $C_{s} = 1/92, 9/190, 9/120$ for $s=1, 0.5, 0.25$ respectively and $W_s$]{{\includegraphics[width=10cm]{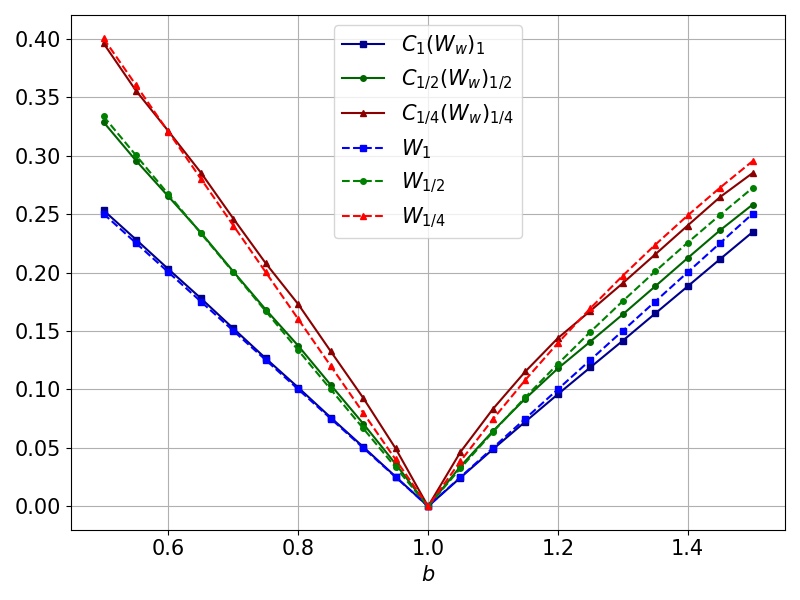} }}%
\caption{Normalized numerical results of three wavelet $s$-Wasserstein distances in comparison with the classical $s$-Wasserstein distance for dilations of the uniform probability measure on $[1,2]$}
\end{figure}
Both \hyperref[Fig2(a)]{Fig. 2 (a)} and \hyperref[Fig2(b)]{Fig. 2 (b)} show some abnormal behavior of the original and reformulated $s$-Wasserstein distances between the uniform probability measure and its highly contracted or stretched counterparts. On the other hand, \hyperref[Fig2(c)]{Fig. 2 (c)} again provides the best performance, capturing well the behavior of $W_s$ under dilations.
\subsubsection{Translations and dilations of the bump probability measure}

\noindent The next two simulations below are analogous to the two above, but with the indicator function replaced by a smooth bump function. The intention is to show the behavior of the wavelet distances when the mass of the measure is not evenly distributed over space and when its density is smooth.\\ 

\noindent \textbf{Simulation three}. (Translations of the bump probability measure on $[0,1]$) \phantomsection\label{Simulation three} Let \begin{equation}\tilde{p}(x) = 
\begin{cases}
C\exp\left( -\frac{1}{1 - (2x)^2}\right), & \text{ if } x \in (-1,1), \\
0, & \text{ if } x\in \mathbb{R}\setminus \{(-1,1)\}, 
\end{cases}
\label{eq4.2}\tag{eq4.2}
\end{equation}
be the bump function supported on $[-1/2,1/2]$, where $C>0$ is chosen so it has unit $L^1$ norm. Let $p(x)\vcentcolon = \tilde{p}(x-1/2)$ and let $p_{a}(x)\vcentcolon = p(x-a)$ be the translation of $p$ by $a>0$. Three wavelet $s$-Wasserstein distances between $p$ and $p_{a}$ for 20 evenly spaces values of $a\in[0,2]$ are computed for three different $s$ and compared to the classical $s$-Wasserstein distance. The $j_0$ and $M$ are chosen to be $-11$ and $22$.
\begin{figure}[H]
\phantomsection\label{Fig3(a)}
\centering
\subfloat[\centering $C_{s}^{\text{o}}(W_{\text{w,o}})_s$ with $C_{s}^{\text{o}} = 1/5, 2/11, 1/6$ for $s=1, 0.5, 0.25$ respectively and $W_s$]{{\includegraphics[width=10cm]{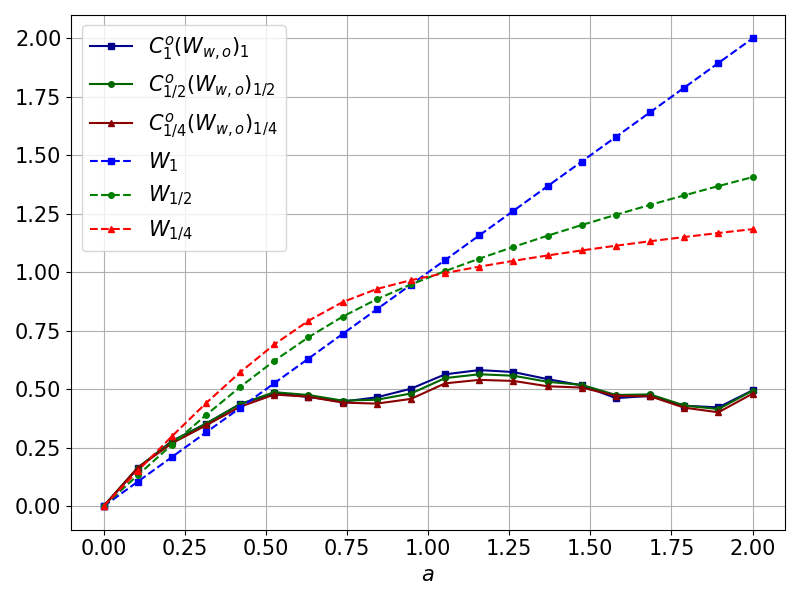} }}%
\end{figure}
\begin{figure}[H]
\phantomsection\label{Fig3(b)}
\centering
\subfloat[\centering $C_{s}^{\text{a}}(W_{\text{w,a}})_s$ with $C_{s}^{\text{a}} = 1/10, 1/7, 4/25$ for $s=1, 0.5, 0.25$ respectively and $W_s$]{{\includegraphics[width=10cm]{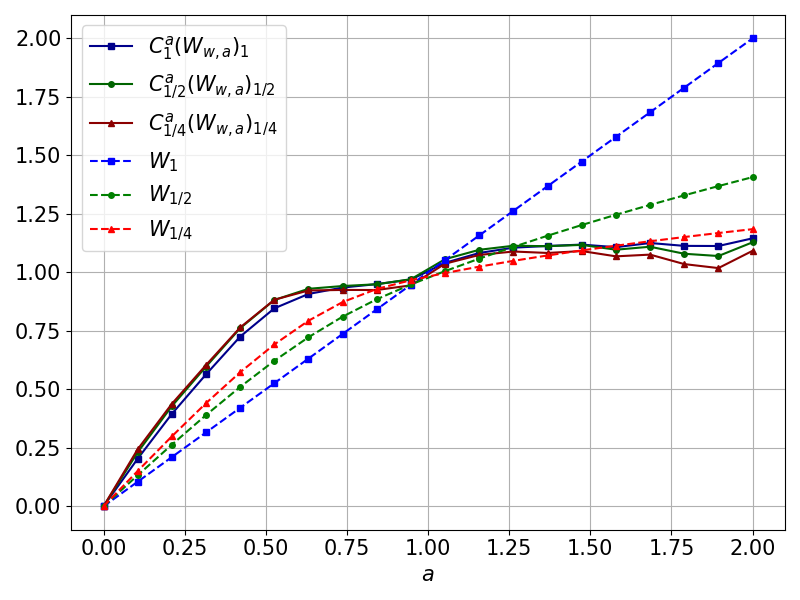} }}%
\end{figure}
\begin{figure}[H]
\phantomsection\label{Fig3(c)}
\centering
\subfloat[\centering $C_{s}(W_{\text{w}})_s$ with $C_{s} = 1/115, 1/25, 1/15$ for $s=1, 0.5, 0.25$ respectively and $W_s$]{{\includegraphics[width=10cm]{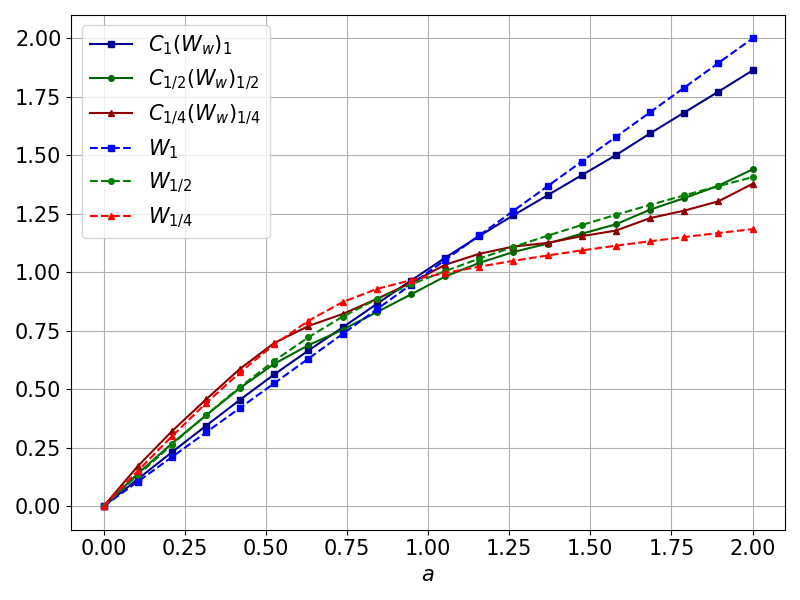} }}%
    \caption{Normalized numerical results of three wavelet $s$-Wasserstein distances in comparison with the classical $s$-Wasserstein distance for translations of the bump probability measure on $[0,1]$}
    \label{fig:example}%
\end{figure}
Similar to the case in Simulation \hyperref[Simulation one]{one}, we observe that the new wavelet s-Wasserstein distance, shown in Figure \hyperref[Fig3(c)]{Fig. 3(c)}, outperforms both the original and reformulated wavelet s-Wasserstein distances in \hyperref[Fig3(a)]{Fig. 3(a)} and \hyperref[Fig3(b)]{Fig. 3(b)} in capturing the behavior of $W_s$ under translations.\\

\noindent \textbf{Simulation four}. (Dilations of the bump probability measure on $[1,2]$) \phantomsection\label{Simulation four}Let $p(x)$ be $p'(x-3/2)$, where $\tilde{p}(x)$ is the function defined in (\ref{eq4.2}). Then let $p_b(x) \vcentcolon= p'((x-3/2)/b)$ be the dilation of $p$ around its center $3/2$ by $b>0$. Three wavelets $s$-Wasserstein distances between $p$ and $p_{b}$ for 20 evenly spaces values of $b\in[1/2,3/2]$ are computed for three different $s$ with the comparison to the classical $s$-Wasserstein distance. The $j_0$ and $M$ are selected to be $-9$ and $22$.
\begin{figure}[H]
\phantomsection\label{Fig4(a)}
\centering
\subfloat[\centering $C_{s}^{\text{o}}(W_{\text{w,o}})_s$ with $C_{s}^{\text{o}} = 1/13, 3/26, 3/22$ for $s=1, 0.5, 0.25$ respectively and $W_s$]{{\includegraphics[width=10cm]{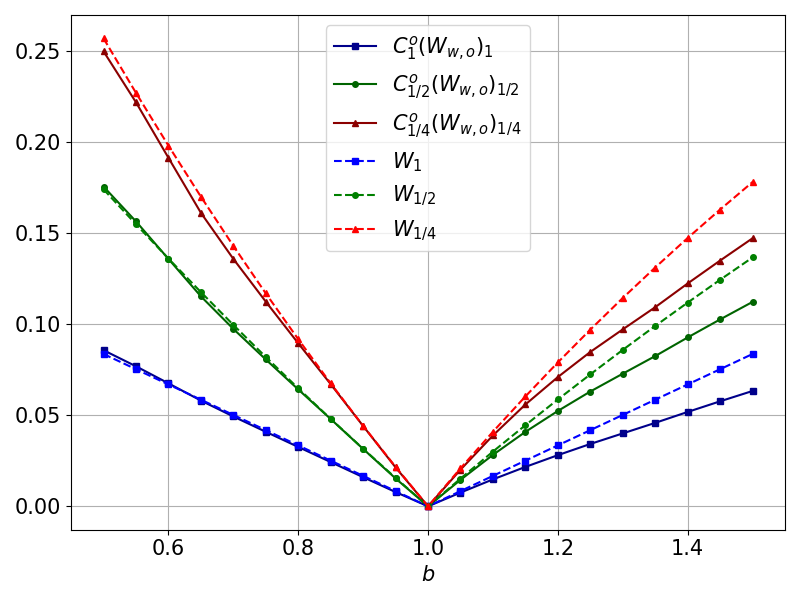} }}%
\end{figure}
\begin{figure}[H]
\phantomsection\label{Fig4(b)}
\centering
\subfloat[\centering $C_{s}^{\text{a}}(W_{\text{w,o}})_s$ with $C_{s}^{\text{a}} = 9/160, 1/10, 4/31$ for $s=1, 0.5, 0.25$ respectively and $W_s$]{{\includegraphics[width=10cm]{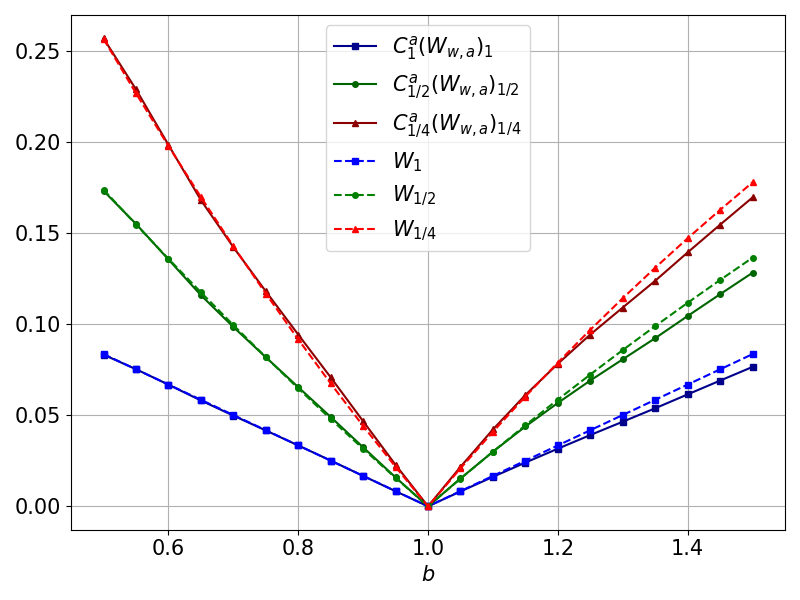} }}%
\end{figure}
\begin{figure}[H]
\phantomsection\label{Fig4(c)}
\centering
\subfloat[\centering $C_{s}(W_{\text{w}})_s$ with $C_{s} = 1/22, 5/54, 5/41$ for $s=1, 0.5, 0.25$ respectively and $W_s$]{{\includegraphics[width=10cm]{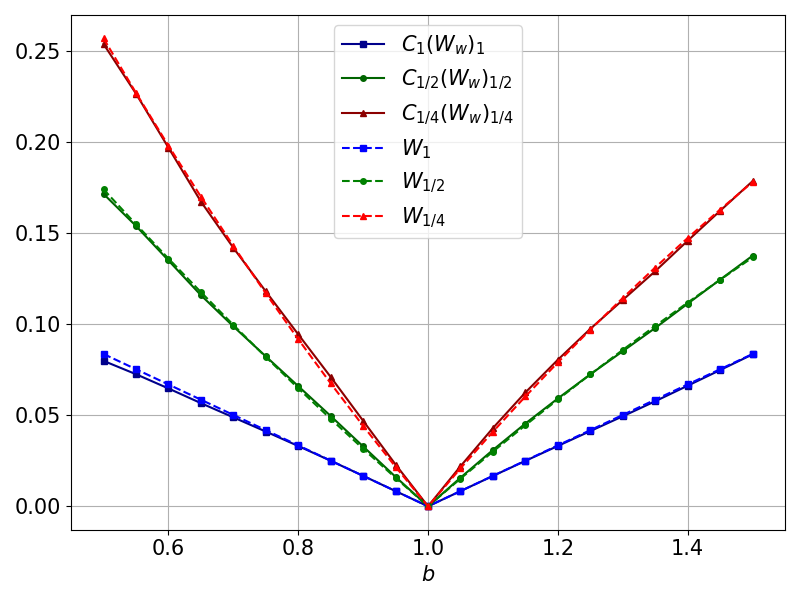} }}%
\caption{Normalized numerical results of three wavelet $s$-Wasserstein distances in comparison with the classical $s$-Wasserstein distance for dilations of the bump probability measure on $[1,2]$}
\end{figure}
Three subfigures, \hyperref[Fig4(a)]{Fig. 4 (a)}, \hyperref[Fig4(b)]{Fig. 4 (b)}, and \hyperref[Fig4(c)]{Fig. 4 (c)}, show that for dilations of smooth bump probability measures, all three formulations of wavelet $s$-Wasserstein distances exhibit reasonable performance in capturing the behavior of the classical $s$-Wasserstein distance. However, still, the new wavelet $s$-Wasserstein distance is the best: it displays exact coincidence after normalization.

\subsection{Numerical dependencies of new wavelet $s$-Wasserstein distance on choice of wavelet system and other parameters}\label{subsec4.3}
In this subsection, we use the same collections of functions as in the previous simulations, but vary other variables to give a better sense of the practical performance of $(W_{\text{wav}}^{(\psi,j_0)})_s$. 

\subsubsection{Dependencies on wavelet systems}
\paragraph{Regular wavelet systems}
First, while the example in Remark \hyperref[Remark 3.19]{3.19} essentially demonstrates that under different wavelet systems, the wavelet $s$-Wasserstein distance can be quite different between certain artificially constructed measures targeting the systems, practitioners might be interested in how $(W_{\text{wav}}^{(\psi,j_0)})_s$ differs in ordinary ``simple'' collections of functions, after normalization. Thus, we redo Simulations \hyperref[Simulation three]{three} and \hyperref[Simulation four]{four} for $(W_{\text{wav,app}}^{(\psi,j_0,M)})_s$, but using different wavelet systems satisfying the conditions in Theorem \hyperref[Theorem 3.13]{3.13}.\\

\noindent \textbf{Simulation five}. (Translations and dilations of the bump probability measure for $(W_{\text{wav,app}}^{(\psi,j_0,M)})_s$ in different wavelet systems) The initial data sets are precisely the ones in Simulations \hyperref[Simulation three]{three} and \hyperref[Simulation four]{four}, and so are the choices of $j_0$ and $M$.\footnote{Technically, $j_0$ and $M$ should be selected after one fixes a wavelet system. However, for comparison purpose, we fix them beforehand in this particular case.} Three wavelet systems, Daubechies 4 wavelet (\texttt{db4}), Daubechies 10 wavelet (\texttt{db10}), and Daubechies 20 (\texttt{db20}), are selected;\footnote{This selection is not entirely random. \texttt{db4} is the Daubechies wavelet with the shortest filter bank satisfying the smoothness condition in Theorem \hyperref[Theorem 3.13]{3.13}. \texttt{db20} on the other hand has a quite long one.} the $s$ values are chosen only to be 1 and 1/2, to avoid over-cluttering on the plots. 
\begin{figure}[H]
\phantomsection\label{Fig5(a)}
\centering
\subfloat[\centering Redo of Simulation three for $C_{(\psi,s)}(W_{\text{w}}^{(\psi)})_s$ with $C_{(\text{db4},1)} = 1/73, C_{(\text{db4},1/2)} = 1/15, C_{(\text{db10},1)} = 1/115, C_{(\text{db10},1/2)} = 1/25, 
C_{(\text{db20},1)} = 1/121, 
C_{(\text{db20},1/2)} = 1/30$]{{\includegraphics[width=10cm]{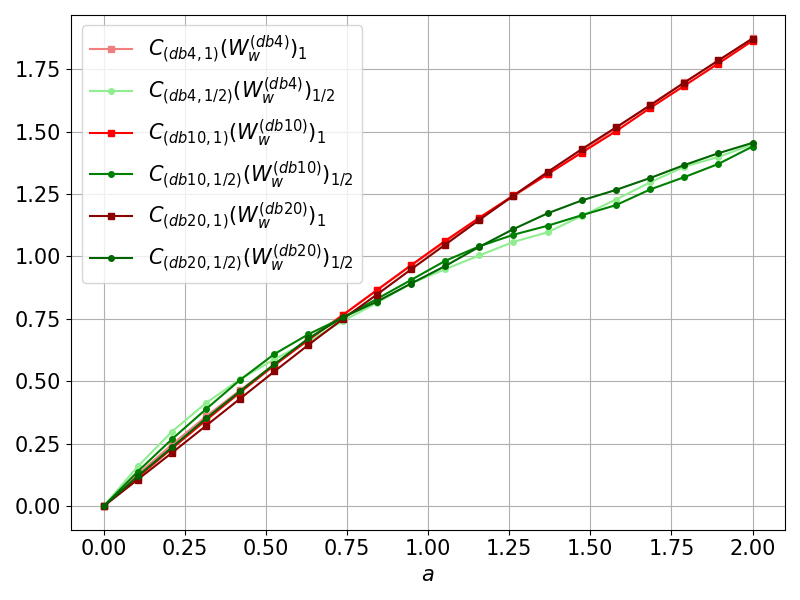} }}%
\end{figure}
\begin{figure}[H]
\phantomsection\label{Fig5(b)}
\centering
\subfloat[\centering Redo of simulation four for $C_{(\psi,s)}(W_{\text{wav}}^{(\psi)})_s$ with $C_{(\text{db4},1)} = 1/25, 
C_{(\text{db4},1/2)} = 5/53, C_{(\text{db10},1)} = 1/22,
C_{(\text{db10},1/2)} = 5/54, C_{(\text{db20},1)} = 1/31,
C_{(\text{db20},1/2)} = 1/15$]{{\includegraphics[width=10cm]{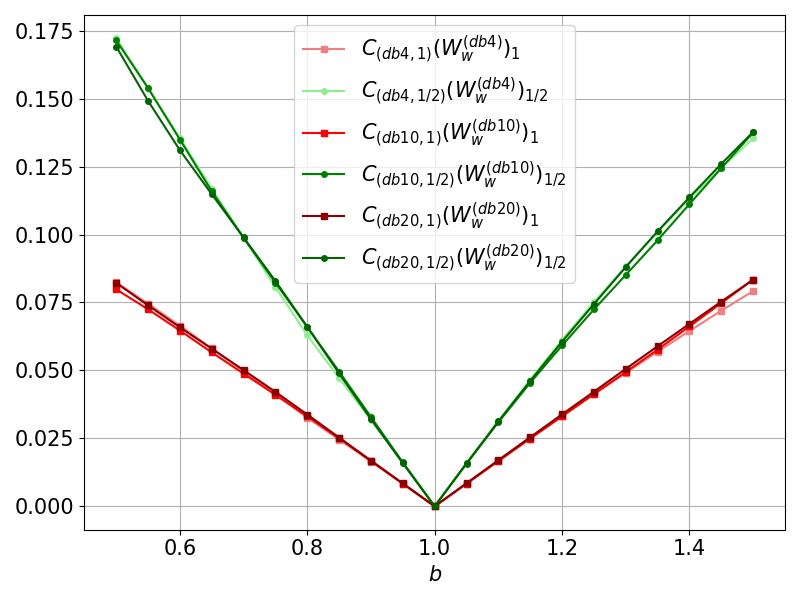} }}%
\caption{Normalized numerical results of the new wavelet $s$-Wasserstein distance under different wavelet systems for translations and dilations of the bump probability measure}
\end{figure}
From \hyperref[Fig5(a)]{Fig. 5 (a)} and \hyperref[Fig5(b)]{Fig. 5 (b)}, it is clear that the choice of wavelet systems here matters little, after normalization.

\paragraph{Nonregular wavelet systems}
One might also be interested in the numerical performances of $(W_{\text{wav}}^{(\psi,j_0)})_s$ for systems that do not meet the smoothness condition specified in Theorem \hyperref[Theorem 3.13]{3.13}. Consequently, we redo Simulations \hyperref[Simulation three]{three} and \hyperref[Simulation four]{four} for $(W_{\text{wav,app}}^{(\psi,j_0,M)})_s$, using wavelet systems that fail to satisfy such condition.\footnote{Additionally, one could explore systems that are not compactly supported; however, in numerical computations, they behave similarly to compactly supported systems, as computers can only perform finite sums.}\\

\noindent \textbf{Simulation six}. (Translations and dilations of the bump measure for $(W_{\text{wav,app}}^{(\psi,j_0,M)})_s$ in different wavelet systems (II)) The initial data sets are precisely the ones in Simulation \hyperref[Simulation three]{three} and \hyperref[Simulation four]{four}, and so are the choices of $j_0$ and $M$. Three wavelet systems, Haar wavelet (\texttt{haar}), Daubechies 2 wavelet (\texttt{db2}), and Daubechies 10 (\texttt{db10}), are selected; the first two systems does not meet the smoothness condition while the last one does; the inclusion of the last one is for comparison. The $s$ values are chosen only on 1 and 1/2 to avoid over-cluttering on the plots. 
\begin{figure}[H]
\phantomsection\label{Fig6(a)}
\centering
\subfloat[\centering Redo of simulation three for $C_{(\psi,s)}(W_{\text{w}}^{(\psi)})_s$ with $C_{(\text{haar},1)} = 1/50, C_{(\text{haar},1/2)} = 1/19, C_{(\text{db2},1)} = 1/130, C_{(\text{db2},1/2)} = 1/22, 
C_{(\text{db10},1)} = 1/115, 
C_{(\text{db10},1/2)} = 1/25$]{{\includegraphics[width=10cm]{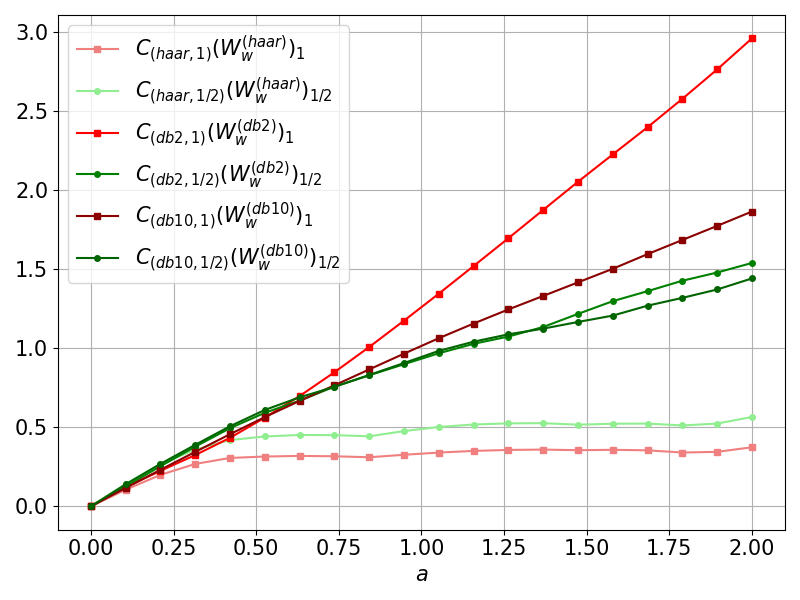} }}%
\end{figure}
\begin{figure}[H]
\phantomsection\label{Fig6(b)}
\centering
\subfloat[\centering Redo of simulation four for $C_{(\psi,s)}(W_{\text{wav}}^{(\psi)})_s$ with $C_{(\text{haar},1)} = 1/5, 
C_{(\text{haar},1/2)} = 3/14, C_{(\text{db2},1)} = 1/45,
C_{(\text{db2},1/2)} = 1/13, C_{(\text{db10},1)} = 1/22,
C_{(\text{db10},1/2)} = 5/54$]{{\includegraphics[width=10cm]{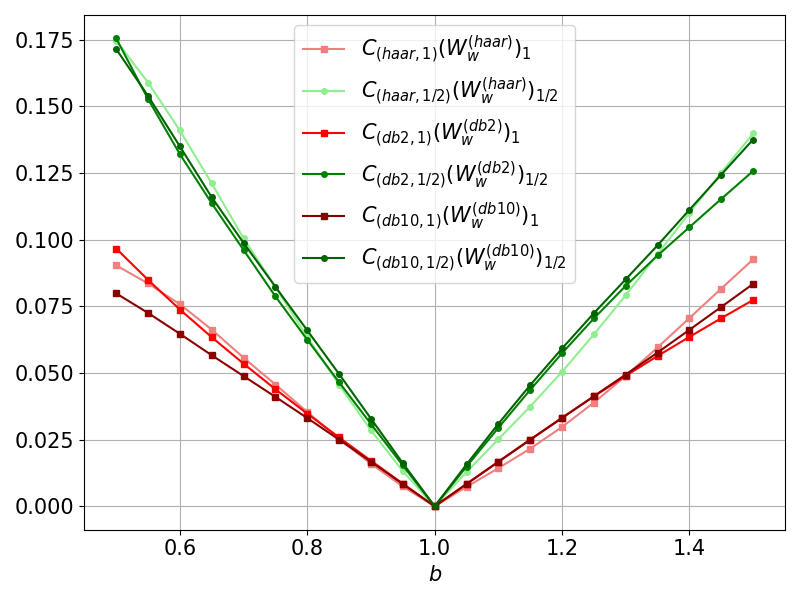} }}%
\caption{Normalized numerical results of the new wavelet $s$-Wasserstein distance under wavelet systems for translations and dilations of the bump probability measure (II)}
\end{figure}
While \hyperref[Fig6(b)]{Fig. 6 (b)} shows that the results are relatively stable for dilations of the bump probability measures under those wavelet systems, \hyperref[Fig6(a)]{Fig. 6 (a)} displays abnormal behavior of the new wavelet $s$-Wasserstein distance for translations under the Haar and db2 wavelet systems, which fail to satisfy the regularity hypotheses of our theorem. This shows that wavelet systems that do not satisfy the smoothness condition, like the popular Haar wavelet system, should be avoided.

\subsubsection{Dependencies on number of coefficients}
\paragraph{Lowest level of coefficients being summed}
The example in Remark \hyperref[Remark 3.19]{3.19} also demonstrates that, under different choices of $j_0$, the wavelet $s$-Wasserstein distance can be quite unpredictable among certain artificially constructed measures that target the system (after fixing the wavelet system). Practitioners might also be interested in observing, after normalization, how $(W_{\text{wav}}^{(\psi,j_0)})_s$ varies across different $j_0$ for ``simple'' data sets. Thus, we redo Simulation \hyperref[Simulation three]{three} and \hyperref[Simulation four]{four} for $(W_{\text{wav,app}}^{(\psi,j_0,M)})_s$, but under different values of $j_0$. Importantly, to isolate the variable of interest --- specifically, the lowest level of all coefficients being summed --- we aim to fix the highest level coefficients summed. This is not merely fixing $M$; rather, $M$ must be adjusted accordingly so that $j_0+M$ is a constant, which can be clearly seen in the steps in Section \ref{subsec4.1}.\\

\noindent \textbf{Simulation seven}\phantomsection\label{Simulation seven}. (Translations and dilations of the bump probability measure for $(W_{\text{wav,app}}^{(\psi,j_0,M)})_s$ with different $j_0$ and fixed $j_0+M$) The initial data sets are precisely the ones in Simulation \hyperref[Simulation three]{three} and \hyperref[Simulation four]{four}, and so is the choice of the wavelet system. For the redo of the Simulation \hyperref[Simulation three]{three}, $j_0$ is selected to be $-5, -8,-11$, with $M$ being $16, 19, 22$ respectively; for the redo of the Simulation \hyperref[Simulation four]{four}, $j_0$ is selected to be $-3, -6,-9$, with $M$ being $16, 19, 22$ respectively. The $s$ values are chosen only to be 1 and 1/2, to avoid over-cluttering on the plots.
\begin{figure}[H]
\phantomsection\label{Fig7(a)}
\centering
\subfloat[\centering Redo of simulation three for $C_{(s,j_0,M)}(W_{\text{w}}^{(j_0,M)})_s$ with $C_{(1,-5,16)} = 1/50, C_{(1/2,-5,16)} = 1/21, C_{(1,-8,19)} = 1/85, C_{(1/2,-8,19)} = 1/24, 
C_{(1,-11,22)} = 1/115, 
C_{(1/2,-11,22)} = 1/25$]{\includegraphics[width=10cm]{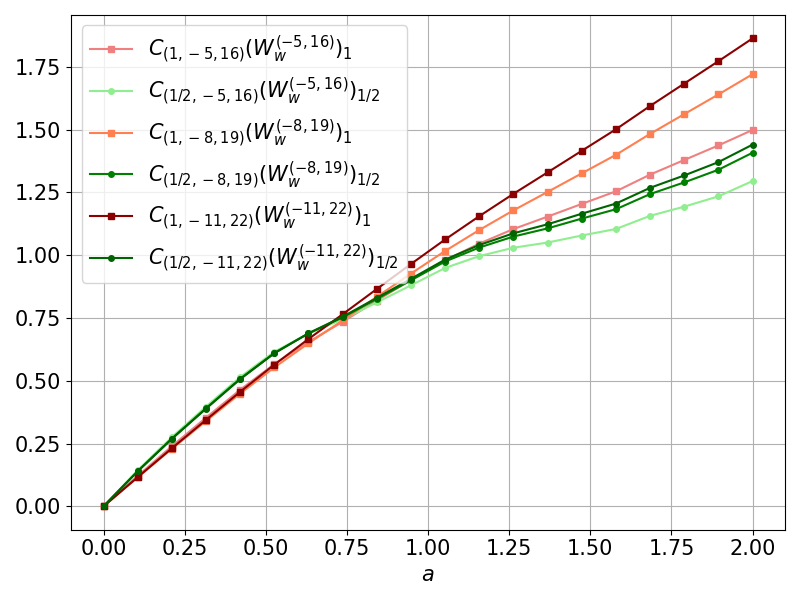}}%
\end{figure}
\vspace{-1.0em}
\begin{figure}[H]
\phantomsection\label{Fig7(b)}
\centering
    \subfloat[\centering Redo of simulation four for $C_{(s,j_0,M)}(W_{\text{w}}^{(j_0,M)})_s$ with $C_{(1,-3,16)} = 1/21, 
    C_{(1/2,-3,16)} = 5/54, C_{(1,-6,19)} = 1/22,
    C_{(1/2,-6,19)} = 5/54, C_{(1,-9,22)} = 1/22,
    C_{(1/2,-9,22)} = 5/54$]{\includegraphics[width=10cm]{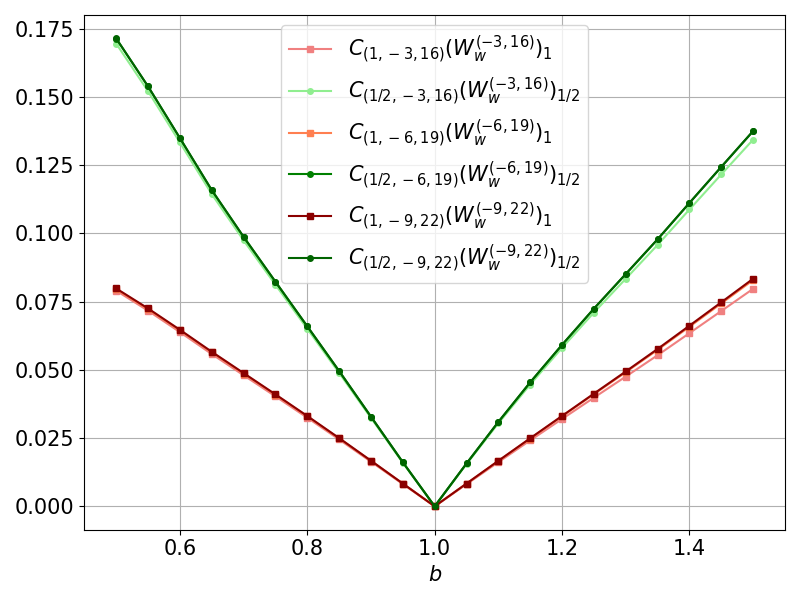}}%
    \caption{Normalized numerical results of the new wavelet $s$-Wasserstein distances under different choices of $j_0$, for translations and dilations of the bump probability measure}
\end{figure}
 From \hyperref[Fig7(b)]{Fig. 7 (b)}, one sees that the results are stable for dilations of the bump probability measures under different choices of $j_0$. On the other hand, \hyperref[Fig7(a)]{Fig. 7 (a)} suggests, the new wavelet $s$-Wasserstein distance is more sensitive to the increase of choice of $j_0$, particularly when the translation is large. More negative values of $j_0$ lead to more accurate approximation of the translation behavior of the true $W_s$ distance.

\paragraph{Highest level of coefficients being summed}
Finally, for completeness, we also redo Simulation \hyperref[Simulation three]{three} and \hyperref[Simulation four]{four} for $(W_{\text{wav,app}}^{(\psi,j_0,M)})_s$, but under different numbers of levels of positive level coefficients. This time, the adjustment involves varying only $M$ while fixing $j_0$. It is important to note that there is no ``symmetry'' here with the previous simulations: as developed in Section \ref{subsec3.4}, the parameter $j_0$ is fundamental to the wavelet Wasserstein $s$-distance; conversely, as seen from Section \ref{subsec4.1}, choosing a large $M$ primarily ensures adequate approximations for both the initialization and final summation. Nevertheless, selecting an appropriate $M$ is crucial from a computational perspective, as it directly determines the level of decomposition in DWT and the number of coefficients generated, which would significantly impact computational time. \\

\noindent \textbf{Simulation eight}. (Translations and dilations of the bump probability measure for $(W_{\text{wav,app}}^{(\psi,j_0,M)})_s$ with different $M$) The initial data sets are precisely the ones in Simulation \hyperref[Simulation three]{three} and \hyperref[Simulation four]{four}, and so are the choices of the wavelet system and $j_0$. For both of the simulations, $M$ is selected to be $16, 19, 22$, respectively. The $s$ values are chosen only to be 1 and 1/2, to avoid over-cluttering on the plots.
\begin{figure}[H]
\phantomsection\label{Fig8(a)}
\centering
\subfloat[\centering Redo of simulation three for $C_{(s,M)}(W_{\text{w}}^{(M)})_s$ with $C_{(1,16)} = 1/115, C_{(1/2,16)} = 1/25, C_{(1,19)} = 1/115, C_{(1/2,19)} = 1/25, 
C_{(1,-22)} = 1/115, 
C_{(1/2,22)} = 1/25$]{{\includegraphics[width=10cm]{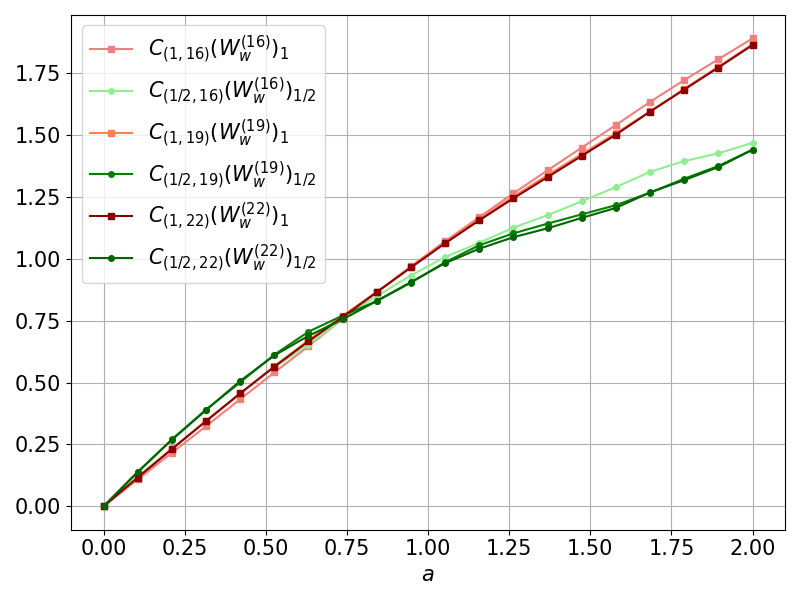} }}%
\end{figure}
\begin{figure}[H]
\phantomsection\label{Fig8(b)}
\centering
\subfloat[\centering Redo of simulation four for $C_{(s,M)}(W_{\text{w}}^{(M)})_s$ with $C_{(1,16)} = 1/24, 
C_{(1/2,16)} = 5/54, C_{(1,19)} = 1/22,
C_{(1/2,19)} = 5/54, C_{(1,22)} = 1/22,
C_{(1/2,22)} = 5/54$]{{\includegraphics[width=10cm]{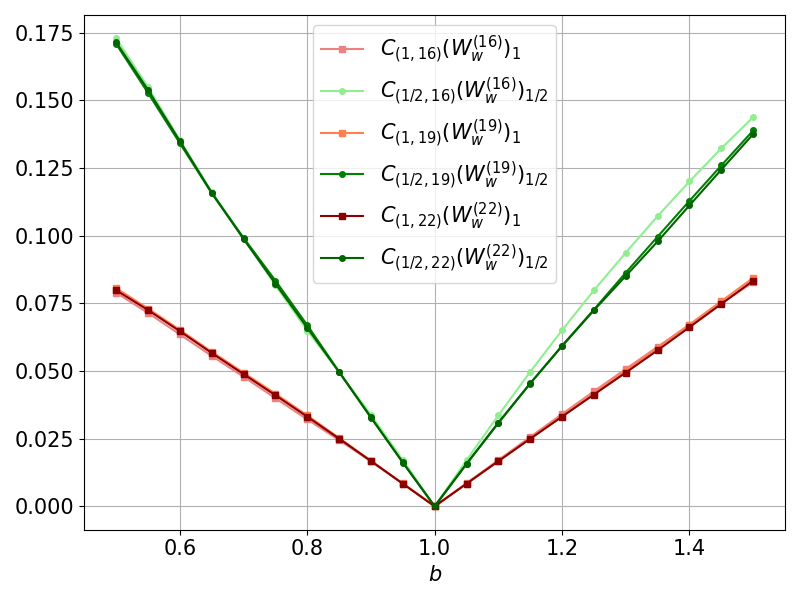} }}%
\caption{Normalized numerical results of the new wavelet $s$-Wasserstein distances under choice of $M$ for translations and dilations of the bump probability measure}
\end{figure}
\hyperref[Fig8(a)]{Fig. 8 (a)} and \hyperref[Fig8(b)]{Fig. 8 (b)} show that the results are relatively stable under different choice of $M$, at least for $M$ not too small. It is worth noting that, across different $M$, our choices of the normalizing constants are almost unchanged, meaning that one may concern less about $M$ when picking the normalizing constants in simulations.\\

\noindent \emph{Acknowledgements.} The authors would like to thank Will Leeb for an interesting discussion that inspired this project. They would also like to thank Caroline Moosmüller for a discussion on the methods for embedding Wasserstein distances into a linear space. 

K. Craig’s work was supported by NSF DMS grant 2145900.
\bibliography{sn-bibliography}

\end{document}